\documentclass[referee]{svjour3}
\usepackage{bm}
\usepackage{graphicx}
\usepackage{amsfonts}
\usepackage{amssymb}
\usepackage{amsmath}
\usepackage{bbm}
\usepackage{latexsym}
\usepackage{hyperref}

\newcommand{\bth}{\begin{theorem}}
\newcommand{\enth}{\end{theorem}}

\newcommand{\bp}{\begin{proposition}}
\newcommand{\ep}{\end{proposition}}

\newtheorem{corrollary}{Corrollary}
\newcommand{\bcor}{\begin{corrollary}}
\newcommand{\ecor}{\end{corrollary}}
\newcommand{\blem}{\begin{lemma}}
\newcommand{\elem}{\end{lemma}}
\newcommand{\bpb}{\begin{question}}
\newcommand{\epb}{\end{question}}

\newcommand{\bdf}{\begin{definition}}
\newcommand{\edf}{\end{definition}}
\newcommand{\bpr} {\begin{proof}}
\newcommand{\epr}{\qed \end{proof}}

\newcommand{\bd}{\begin{document}}
\newcommand{\ed}{\end{document}}
\newcommand{\beq}{\begin{equation}}
\newcommand{\eeq}{\end{equation}}
\newcommand{\bef}{\begin{figure}}
\newcommand{\enf}{\end{figure}}
\newcommand{\bea}{\begin{eqnarray}}
\newcommand{\eea}{\end{eqnarray}}
\newcommand{\baR}{\begin{array}}
\newcommand{\eaR}{\end{array}}
\newcommand{\bc}{\begin{center}}
\newcommand{\ec}{\end{center}}
\newcommand{\ben}{\begin{enumerate}}
\newcommand{\een}{\end{enumerate}}
\newcommand{\bit}{\begin{itemize}}
\newcommand{\eit}{\end{itemize}}
\newcommand{\su}{\section}
\newcommand{\ssu}{\subsection}
\newcommand{\sssu}{\subsubsection}
\newcommand{\nid}{\noindent}
\newcommand{\nnb}{\nonumber}

\newcommand\be{{\bf e}}
\newcommand\bv{{\bf v}}
\newcommand\h{{\bf h}}
\newcommand\bom{\mbox{{\boldmath $\omega$}}}
\newcommand\tom{\omega}
\newcommand\cM{{\cal M}}
\newcommand\cE{{\cal E}}
\newcommand\hM{{\hat{\cal M}}}
\newcommand\Mom{\cM_{\omega}}
\newcommand\on{\seq{\omega}{0}{n}}
\newcommand\onm{\seq{\omega}{0}{n-1}}
\newcommand\onmd{\left[\omega\right]_{n-2}}
\newcommand\ot{\left[\omega\right]_t}
\newcommand\otm{\left[\omega\right]_{t-1}}
\newcommand\Mon{\cM_{\on}}
\newcommand\Monm{\cM_{\onm}}
\newcommand\hMon{\hat{\cM}_{\on}}
\newcommand\hMot{\hat{\cM}_{\ot}}
\newcommand\oMon{\stackrel{\circ}{\cM}_{\on}}
\newcommand\ohMon{\stackrel{\circ}{\hM}_{\on}}
\newcommand\Anion{A_{i,n,\on}}
\newcommand\bAnion{\bar{A}_{i,n,\on}}
\newcommand\Fon{\F_{\on}}
\newcommand\Fot{\F_{\ot}}
\newcommand\hF{\hat{\F}}
\newcommand\hFon{\hF_{\on}}
\newcommand\hFot{\hF_{\ot}}
\newcommand\V{{\bf V}}
\newcommand\sV{\bf{\tiny{V}}}
\newcommand\f{{\bf{f}}}
\newcommand\hV{{\hat{\V}}}
\newcommand\hS{{\hat{\cS}}}
\newcommand\hml{{\hat{\mu}_0}}
\newcommand\piu{\pi^{(u)}}
\newcommand\pis{\pi^{(s)}}
\newcommand\tin{\tau_i(\omega,n)}

\newcommand\cEta{\cC(\bEta)}

\newcommand\bEtap{\mbox{{\boldmath $\eta'$}}}
\newcommand\cEtap{\cC(\bEtap)}
\newcommand\betat{\tilde{\bEta}}
\newcommand\betatp{\tilde{\bEtap}}
\newcommand\betatpl{\tilde{\bEta}^+}
\newcommand\bbeta{\mbox{{\boldmath $\beta$}}}
\newcommand\mub{\mu_{\bbeta}}
\newcommand\mus{\mu^\ast}
\newcommand\Ps{P^\ast}
\newcommand\Fb{\cF_{\bbeta}}
\newcommand\bxi{\bm{\xi}}
\newcommand\balpha{{\bm{\alpha}}}
\newcommand\blambda{{\bm{\lambda}}}
\newcommand\bphi{{\bm{\phi}}}
\newcommand\B{{\bf B}}
\newcommand\Vp{{\bf V'}}
\newcommand\W{{\bf W}}
\newcommand\Z{{\bf Z}}
\newcommand\F{{\bf F}}
\newcommand\Fe{\F_{\bEta}}
\newcommand\Fei{F_{\bEta,i}}
\newcommand\G{{\bf G}}
\newcommand\I{{\bf I}}
\newcommand\bz{{\bf 0}}
\newcommand\bzs{\left\{\bz\right\}}
\renewcommand\S{{\bf S}}
\newcommand\Td{T\left(\dAS\right)}
\newcommand\X{{\bf X}}
\newcommand\x{{\bf x}}
\newcommand\y{{\bf y}}
\newcommand\xt{\X_t(\V)}
\newcommand\cktv{\chi_{k,t}(\V)}
\newcommand\cA{{\cal A}}
\newcommand\cB{{\cal B}}
\newcommand\cH{{\cal H}}
\newcommand\cL{{\cal L}}
\newcommand\cN{{\cal N}}
\newcommand\cO{{\cal O}}
\newcommand\bme{\cB_\cM(\epsilon)}
\newcommand\Bev{\cB(\V,\epsilon)}
\newcommand\coMe{\stackrel{\circ}{\cMe}}
\newcommand\beV{\cB(\V,\epsilon)}
\newcommand\bmd{\cB_\cM(\delta)}
\newcommand\bmdp{\cB_\cM(\delta')}
\newcommand\cC{{\cal C}}
\newcommand\cD{{\cal D}}
\newcommand\cJ{{\cal J}}
\newcommand\cR{{\cal R}}
\newcommand\De{\cD(\bEta)}
\newcommand\bDe{\bar{\cD}(\bEta)}
\newcommand\Dep{\cD(\bEtap)}
\newcommand\bDep{\bar{\cD}(\bEtap)}
\newcommand\cG{{\cal G}}
\newcommand\GW{\cG_{(\cW,\Ie)}}
\newcommand\IW{\cI_{(\cW,\Ie)}}
\newcommand\cF{{\cal F}}
\newcommand\bcF{\bar{\cF}}
\newcommand\oGF{{\omega(\GF)}}
\newcommand\oGFTe{{\omega(\GFTe)}}
\newcommand\oF{{\omega_{\cF}(\GF)}}
\newcommand\oFTe{{\omega_{\cF}(\GFTe)}}
\newcommand\Bie{\cB_{\infty,\epsilon}}
\newcommand\BTe{\cB_{T,\epsilon}}
\newcommand\BTep{\cB_{{T+1},\epsilon}}
\newcommand\BTeu{\cB_{T,\epsilon}^{(1)}}
\newcommand\BTed{\cB_{T,\epsilon}^{(2)}}
\newcommand\DFFTV{\left. D\FTF\right|_{\V}}
\newcommand\FF{{\F_{\cF}}}
\newcommand\FTF{{\FT_{\cF}}}
\newcommand\FFb{{\F_{\bcF}}}
\newcommand\VF{{\V_{\cF}}}
\newcommand\VFb{{\V_{\bcF}}}
\newcommand\ft{\cF_t(\V)}
\newcommand\ftmu{\cF_{t-1}(\V)}
\newcommand\ftmk{\cF_{t-k}(\V)}
\newcommand\ftmkpu{\cF_{t-k+1}(\V)}
\newcommand\bft{\bar{\cF}_t(\V)}
\newcommand\bftmu{\bar{\cF}_{t-1}(\V)}
\newcommand\bftmk{\bar{\cF}_{t-k}(\V)}
\newcommand\bftmkpu{\bar{\cF}_{t-k+1}(\V)}
\newcommand\cMz{{\cal M_{\bf{0}}}}
\newcommand\cMu{{\cal M_{\bf{1}}}}
\newcommand\cMe{{\cM_{\bEta}}}
\newcommand\cMeu{{\cM_{{\bEta}_1}}}
\newcommand\cMet{{\cM_{{\bEta}_t}}}
\newcommand\cMeT{{\cM_{{\bEta}_T}}}
\newcommand\cMep{{\cM_{\bEtap}}}
\newcommand\cU{{\cal U}}
\newcommand\cS{{\cal S}}
\newcommand\cI{{\cal I}}
\newcommand\cP{{\cal P}}
\newcommand\cQ{{\cal Q}}
\newcommand\cV{{\cal V}}
\newcommand\cW{{\cal W}}
\newcommand\Vpi{V^+_i}
\newcommand\vp{\cV^+}

\newcommand\Is{\I^s}
\newcommand\IsP{\I^{s\ast}}
\newcommand\Ise{\I^s(\bEta)}
\newcommand\IE{\I(\bEta)}
\newcommand\Isv{\I^s(\V)}
\newcommand\Isi{I^s_i}
\newcommand\Isj{I^s_j}
\newcommand\Iie{I_i(\bEta)}
\newcommand\Iiep{I_i(\bEta')}
\newcommand\Ije{I_j(\bEta)}
\newcommand\Ijep{I_j(\bEta')}
\newcommand\Isie{\Isi(\bEta)}
\newcommand\Isje{\Isj(\bEta)}
\newcommand\Isjep{\Isj(\bEtap)}
\newcommand\Isiv{\Isi(\V)}
\newcommand\Isvt{I^s_i(\V(t))}
\newcommand\Ie{{\bf I^{ext}}}
\newcommand\Iei{I_i}
\newcommand\Iej{I^{ext}_j}
\newcommand\tk{\tau^{(k)}}
\newcommand\tkp{\tau^{(k+1)}}
\newcommand\toi{\tau^{(1)}_i(\V)}
\newcommand\toik{\tau^{(k)}_i(\V)}
\newcommand\toike{\tau^{(k)}_i(\betat)}
\newcommand\toikpe{\tau^{(k+1)}_i(\V)}
\newcommand\toikep{\tau^{(k+1)}_i(\betat)}
\newcommand\dik{\delta_i^{(k)}(\V)}
\newcommand\PrW{\cP_\cW}
\newcommand\PrI{\cP_\Ie}
\newcommand\PrWi{\cP_{\cW,\Ie}}
\newcommand\PW{\cP_{(\cW,\Ie)}}
\newcommand\RW{\cR_{(\cW,\Ie)}}
\newcommand\cPT{\cP^{(T)}}
\newcommand\cMT{\cM_{{\bEta_0} \dots {\bEta_T}}}
\newcommand\FT{\F^T}
\newcommand\Ft{\F^t}
\newcommand\Ftp{\F^{t+1}}
\newcommand\Fttk{\F^{t_k}}
\newcommand\Fmu{\F^{-1}}
\newcommand\Fmd{\F^{-2}}
\newcommand\FmT{\F^{-T}}
\newcommand\dAS{d(\oM,\cS)}
\newcommand\dSV{d(\tVp,\cS)}
\newcommand\trV{\left\{\V(t)\right\}_{t \geq 0}}
\newcommand\tV{\tilde{\V}}
\newcommand\tVp{\tilde{\V}^+}
\newcommand\tpi{t^{pre}_i}
\newcommand\tpj{t^{post}_j}
\newcommand\mW{m_{\cW}}
\newcommand\sW{\sigma_{\cW}}
\newcommand\mI{m_{\Ie}}
\newcommand\sI{\sigma_{\Ie}}
\newcommand\sw{\sigma_{(\cW,\Ie)}}
\newcommand\sB{\sigma_B^2}
\newcommand\sd{\sigma^2}
\newcommand\Cst{[\betat]_{s,t}}
\newcommand\defCst{\left\{\betat | \bEta_s=\as, \dots, \bEta_t = \at    \right\}}
\newcommand\as{\mbox{{\boldmath{$\alpha$}}}_s}
\newcommand\at{\mbox{{\boldmath{$\alpha$}}}_t}
\newcommand\moytpsi{\bar{\psi}_\cW}
\newcommand\sn{\left\{1 \dots N \right\}}
\newcommand\GF{\Gamma_\cF}
\newcommand\GFTe{\Gamma_{\cF,T,\epsilon}}
\newcommand\GFTep{\Gamma_{\cF,T,\epsilon'}}
\newcommand\GFTpe{\Gamma_{\cF,{T+1},\epsilon}}
\newcommand\PF{\Pi_{\cF}}
\newcommand\PbF{\Pi_{\bcF}}
\newcommand\oM{\Omega}

\newcommand\wsg{\cW^s(\V)}
\newcommand\wsl{\cW^s_{loc}(\V)}
\newcommand\wslf{\cW^s_{loc}(\F(\V))}
\newcommand\wse{\cW^s_{\epsilon}(\V)}
\newcommand\wset{\cW^s_{\epsilon}(\V(t))}
\newcommand\tjk{\tau_j^{(k)}}
\newcommand\tik{\tau_i^{(k)}}
\newcommand\tiu{\tau_i^{(1)}}
\newcommand\tikp{\tau_i^{(k+1)}}
\newcommand\Vr{V_{reset}}
\newcommand\Vm{V_{min}}
\newcommand\VM{V_{max}}
\newcommand\Vs{V^{+}}
\newcommand\Vsi{V_i^{+}}
\newcommand\SLp{\Sigma_\Lambda^+}
\newcommand\SL{\Sigma_\Lambda}
\newcommand\SW{\Sigma_{(\cW,\Ie)}}
\newcommand\SWp{\Sigma_{(\cW,\Ie)}^+}
\newcommand{\D}{\displaystyle}
\newcommand{\deq}{\stackrel {\rm def}{=}}
\newcommand{\peq}{\stackrel {\rm \mu p.s.}{=}}

\newcommand\cyl[3]{\left[{#1}\right]_{#2}^{#3}}
\newcommand\seq[3]{{#1}_{#2}^{#3}}
\newcommand\prodset[3]{{#1}_{#2}^{#3}}
\newcommand\tinm{\tau_i(\seq{\omega}{-\infty}{n-1})}
\newcommand\uom{\underline\omega}
\newcommand\uX{\underline X}
\newcommand\ucF{\underline \cF}
\newcommand\uomp{\underline\omega'}
\newcommand\tiz{\tau_i(\uom)}
\newcommand{\setZ}{\mathbbm{Z}}
\newcommand{\setN}{\mathbbm{N}}
\newcommand{\setR}{\mathbbm{R}}
\newcommand{\Exp}[1]{\mathbb{E}\left [ #1 \right]}
\newcommand{\ExpB}[1]{\mathbb{E}_B\left [ #1 \right]}
\newcommand{\CovB}[1]{Cov\left [ #1 \right]}
\newcommand\pTo{{\pi_{\omega}^{(T)}}}
\newcommand\tjn{t_j^{(n)}}
\newcommand\tot{\seq{\omega}{-\infty}{t}}
\newcommand\Jito{J_i(t,\tot)}

\begin{document}

\title{A discrete time neural network model with spiking neurons}
\subtitle{II. Dynamics with noise.}
\titlerunning{Discrete time spiking neurons with noise.}

\author{
B. Cessac
}

\institute{
\at Equipe syst\`emes dynamiques, interactions en physique, biologie, chimie,
Laboratoire Jean-Alexandre Dieudonn\'e, Universit\'e de Nice, Parc Valrose, 06000 Nice, France.
\and NeuroMathComp, INRIA, 2004 Route des Lucioles, 06902 Sophia-Antipolis, France.
}


\maketitle

\begin{abstract}
We provide rigorous and exact results characterizing the statistics of 
spike trains  in a network of leaky Integrate-and-Fire neurons,
where time is discrete and where neurons are submitted to noise,
without restriction on the synaptic weights.  We show the
existence and uniqueness of an invariant measure of Gibbs type and discuss its properties.
We also discuss Markovian approximations and relate them to the approaches
currently used in computational neuroscience to analyse experimental spike trains statistics.
 \end{abstract}

\su{Introduction}

The neuronal activity is manifested by the emission of 
\textit{action potentials} or \textit{spikes}. While the shape
of an action potential is essentially constant for a given
neuron, the succession of spikes (\textit{spike train}) that
a neuron is able to emit, depending on its state and in response
to  excitations coming from other neurons or external stimuli,
is simply overwhelming. About twenty different spike trains forms
are classified in the literature \cite{izhikevich:04}.
It is widely believed by the neuroscience community that
spike trains emitted by a neuron assembly constitute
somehow a ``code''  and deciphering this code is
a big challenge \cite{rieke-etal:96}.

Spike train are usually not exactly reproducible
when repeating the same experiment\footnote{Although, retinal responses to a natural image,
seem to be almost reproducible spike by spike \cite{perrinet:08}},
even with a very good control ensuring that the experimental
conditions have not changed. Therefore,
 researchers are seeking statistical
regularities in spike trains. For this, they
 define statistical
indicators such as firing rate, probability of spike
coincidence, spike response function, spike correlations
 (see \cite{rieke-etal:96,dayan-abbott:01,gerstner-kistler:02}
for a comprehensive introduction to spike train analysis).
An early step for ''reading the code'' is therefore to provide
an accurate model for spike train statistics, i.e. a probability
distribution ``fitting at best'' the experimental data and/or matching what neuroscientists believe relevant in neurons communication
via spikes.

For example, it has been long believed that firing rates (the probability that a neuron emits a spike in a certain time interval) were carrying most of the ``information'' exchanged by neurons. As a consequence the canonical statistical model, namely the probability distribution which reproduces the firing rates without additional 
assumptions, is a Bernoulli distribution (possibly with time dependent probabilities), and the probability that a given number of spikes is emitted within a definite time interval is  
 Poisson. Actually, there are many mechanisms in the nervous
system, such as the muscle commands \cite{adrian-zotterman:1926},  working essentially with rates. But more recent experiments evidenced
the role of spikes timing, spike ordering, spike synchronization, 
 in processes such as vision \cite{thorpe:90,thorpe-fize-etal:96,rullen-thorpe:01} or interactions between perception and motion \cite{grammont-riehle:99,riehle-etal:00,grammont-riehle:03}.
Here, one has to consider more elaborated statistical models, 
such as the (weak) pairwise interactions model proposed by Schneidman and collaborators \cite{schneidman-etal:06} in experiments on the salamander retina.
Consequently, there is an intensive activity and wide debates,
focusing on the determination of statistical models of spike train statistics, with a clear evidence: distinct 
statistical models
lead to fundamentally distinct characterizations of the mechanisms
at work in the piece of nervous system under study \cite{nirenberg-latham:03}.  \\

Clearly, obtaining  statistical models from experimental data
or selecting a model among many others 
are difficult task. Forgetting about all the experimental difficulties to obtain
``clean'' data with a good control on parameters experiments, one has
still to solve delicate questions such as the control of finite sampling effects (finite duration, 
finite numbers of experiments),
extrapolation of the probability
distribution characterizing small neural assemblies to a large population of neurons  \cite{roudy-nirenberg-etal:09}, non stationarity, effects of synaptic plasticity or adaptation mechanisms \cite{toyoizumi-etal:07}.
 As a consequence, there is no general recipe to extract a statistical model from data and several approaches
have been proposed \cite{schneidman-etal:06,pouzat-chaffiol:09,marre-boustani-etal:09}.

It appears simpler to characterize spike trains statistics in 
neural networks \textit{models} where one controls exactly
the neural network parameters, the number of involved
neurons, the number of samples, and the duration of the 
experiment (with a possible mathematical extrapolation
to infinite time). Especially, producing analytical
(and when possible, rigorous) results on those statistics
provides clues toward resolving the delicate experimental questions
raised above, with possible outcomes
toward new algorithms for data treatments 
\cite{vasquez-cessac-etal:10}. Obviously, for this, one needs models which are a good compromise
between analytical tractability and biological realism.

Generalized Integrate-and-Fire (gIF) models \cite{rudolph-destexhe:06}
 constitute a good example of this. Besides the fact
that they capture the conductance-based mechanisms for spike generation,
 without
focusing to much on the biological machinery, it has been shown
by authors like Jolivet et al. \cite{jolivet-lewis-etal:04,jolivet-rauch-etal:06}
that they are ``good enough'' to reproduce spike trains
from real neurons. Moreover, these models allow the analytical
characterization of their dynamics \cite{cessac-vieville:08}.
Further simplifications of gIF models lead to the Leaky Integrate-and-Fire (LIF) model, which was in fact the first proposed to model
neuron dynamics (in 1907!) \cite{lapicque:07}. In this setting, prominent mathematical results
on spike statistics in the presence of noise, have been published.
For example, Brunel and Hakim obtained a complete characterization 
of LIF models with noise and strong dilution of synaptic weights,
 using a mean-field approximation and assuming
that the synaptic weights are inhibitory \cite{brunel-hakim:99}.
Also,
Touboul and Faugeras \cite{touboul-faugeras:07} obtained rigorous results
on the probability distribution of inter-spike intervals for one LIF neuron
submitted to noise. They recently extended their results in \cite{touboul-faugeras:09} to networks of IF neurons of several type considering different types of interactions and conclude that the spikes times can be modelled as a Markov chain.\\

In this paper, we proceed along similar lines, although using different methods
and raising 
different conclusions,
and make a complete characterization of
spike train statistics for the discrete-time leaky Integrate-and-Fire  model with  noise and time-independent stimuli.
This is somehow a continuation of the paper \cite{cessac:08}
proposing a complete classification of dynamics in this model, without noise, and of
\cite{cessac-vieville:08} extending these results
to gIF models. It also rigorously supports
the main assumption made in \cite{cessac-rostro-etal:09}
where it was proposed to characterize spike train statistics
in neural networks by Gibbs distributions, with an emphasis on
synaptic plasticity effects.
Here, we propose a framework allowing
to handle dynamics with noise, with possible extensions
to more realistic neural networks models such as gIF (see the 
conclusion section). We emphasize that we do not use any simplifying
assumption in the model. Especially, our results hold for finite-sized
networks, without restriction on the synaptic weights (except that they
are finite)  and all type of
synaptic graph structures are allowed. Also, we are not constrained
by an ad hoc choice of the initial conditions distribution of membrane potentials;
 instead we propose a procedure where this distribution is selected by dynamics
and  is uniquely determined, as we show.

Moreover,
this work attempts to bridge a gap between the mathematical characterization
of spike train statistics and empirical methods or algorithms
currently used by the neuroscience community.
As a consequence, this paper addresses to 2 distinct communities.
On one hand, to specialists
from mathematical statistical physics and ergodic theory, as
far as the mathematics of this paper are concerned. From this point of
view the results exposed here are direct applications of classical results
in ergodic theory.
But, to the best of our knowledge, it is the first time that they are used
in this context. On the other hand, this paper addresses to neuroscientists. Studying a fairly
simple model, from the biological point of view, we nevertheless obtain conclusions
 which could be useful for  the characterization of spike trains in real experiments, 
with concrete applications toward implementation of software for spike train analysis. \\ 

The paper is organized as follows. In section
\ref{Sdef} we define the model and infer some preliminary
results. Especially we compute explicitly the probability
that neurons fire at time $t$ given the past. This defines
a transition probability  which is the main object of our study.
A salient result obtained in section \ref{Sdef} is the fact
that this transition probability is \textit{non Markovian},
it depends on an unbounded past. 
This defines a stochastic
process, known under the name of ``chain with complete connections''
(\cite{maillard:07} and references therein),
which is studied in section \ref{SGibbs}. Especially, we show
that there is a unique invariant probability measure (equilibrium state) whatever
the model-parameters values, which satisfies a variational
principle and is a Gibbs distribution. We also show that the entropy
of the discrete-time leaky Integrate-and-Fire  model with noise is always positive.
In section \ref{SFiniteRange} we propose a Markovian
approximation where memory depends on $R$ time steps. This
approximation allows the 
 computation of the main quantities
used in neuroscience for the characterization of raster
plots statistics. 
The computation of these quantities
is done in section \ref{SSpikeStat}. We also show that, in this
approximation, the equilibrium state can also be obtained
via the Jaynes principle of statistical physics (maximizing
the entropy under constraints), and we discuss in which sense the 
statistical models used in neuroscience community are
approximations whose degree of accuracy can be controlled.

\su{Definitions and preliminary results}\label{Sdef}

\ssu{Model definition.} 

\sssu{The neural network.}  Fix $N > 0$ a positive integer called ``the dimension of the neural network'' (the number
of neurons). 
Let $\cW$ be a $N \times N$ matrix, called ``the matrix of synaptic weights'', 
with entries $W_{ij}$.  It defines
an oriented and signed graph, called ``the neural network associated to $\cW$'',  with  vertices $i=1 \dots N$ called the ``neurons\footnote{Therefore, ``neurons'' are points here, i.e. they have no structure.}''.
There is an oriented edge $j \to i$ whenever $W_{ij} \neq 0$.  $W_{ij}$ is called ``the synaptic weight\footnote{On biological grounds, this corresponds to the maximal amplitude of the post-synaptic potential generated,
at the dendrite connecting the pre-synaptic neuron $j$ to the post-synaptic  neuron $i$, 
when  neuron $j$ emits an action potential.} from neuron $j$ to neuron $i$''. 
The synaptic weight is called ``excitatory'' if $W_{ij}>0$ and ``inhibitory'' if $W_{ij}<0$.
 We assume that the synaptic weights are bounded i.e. 
$W_{ij} \in [W_{min},W_{max}], \forall i,j$ where $-\infty < W_{min} \leq W_{max} < +  \infty$. Moreover, in this paper, the $W_{ij}$'s do not evolve in time.

\sssu{Membrane potential.} Each vertex (neuron) $i$ is characterized by a real variable $V_i$ called the ``membrane potential of neuron $i$''.
Fix a positive real number $\theta>0$ called the ``firing threshold''. Let $Z$ be the function
$Z(x)=\chi(x \geq \theta)$ where $\chi$ is the indicatrix function. Namely,
$Z(x)=1$ whenever $x \geq \theta$ and $Z(x)=0$ otherwise. $Z(V_i)$ is called the ``firing state of neuron $i$''.
When $Z(V_i)=1$ one says that neuron $i$  ``fires'' or ``spikes'' and when $Z(V_i)=0$ neuron $i$
is ``quiescent''. We extend the definition of $Z$ to vectors: $Z(V)$ 
is the vector
with components $Z(V_i), \ i=1 \dots N$.

\sssu{Dynamics.} Fix $\gamma \in [0,1[$,
called the ``leak rate''. The discrete time and synchronous dynamics of our model is given
by:
\beq \label{DNN}
V(t+1)=F(V(t))+ \sigma_B B(t),
\eeq
\nid where $\sigma_B > 0$, $V=\left(V_i\right)_{i=1}^N$ 
is the vector of membrane potentials  and
$F=\left(F_i\right)_{i=1}^N$ with:
%
$$F_i(\V)=\gamma V_i \left(1 - Z[V_i] \right)+ \sum_{j=1}^N W_{ij}Z[V_j]+ \Iei; \qquad i=1 \dots N.$$
%
We assume that initial conditions belong to some compact set in $ \setR^N$ (i.e. the initial
membrane potentials are bounded).
The variable $\Iei$ is called ``an external input applied to neuron $i$''.
We  assume in  this paper that it does not depend on time.

\sssu{Noise.} The vector  $B(t)=\left(B_i(t)\right)_{i=1}^N$ is an additive noise\footnote{On phenomenological grounds, it mimics effects such as noise in synaptic transmission (neurotransmitters diffusion), randomness in ionic channels transitions, or effects of hidden degrees of freedom.}.
 It
has Gaussian identically distributed and
independent entries
$B_i(t)$ with zero mean and variance $1$.
We note $\cN(0,1)$ the standard Gaussian law and :
%
$$\pi(x)=\frac{1}{\sqrt{2\pi}}\int_x^{+\infty} e^{-\frac{u^2}{2}}du.$$
%
The parameter $\sigma_B$ in (\ref{DNN}) tunes the noise amplitude.

\sssu{Interpretation.} 

To the best of our knowledge this model has been first introduced by G. Beslon, O. Mazet
and H. Soula \cite{soula-beslon-etal:06}. 
It belongs to the family of the so-called leaky-Integrate-and-Fire   models
\cite{gerstner-kistler:02}. Its interpretation is the following. 
A neuron ``fires'' i.e. emits an   action potential (or ``spike'')
whenever its membrane potential exceeds the threshold $\theta$. Here a spike is modelled by the function $Z$.
For an isolated neuron, firing corresponds, in the model, to the reset of the membrane potential to a rest value
$V_{rest}=0$.
In a network, each neuron $i$ receives spikes from pre-synaptic neurons.  When a pre-synaptic neuron $j$ emits
a spike this modifies the membrane potential of neuron $i$ by an amount $W_{ij}$. 
Thus, according to eq. (\ref{DNN}), when a neuron fires, it immediately receives inputs from other neurons and from the environment (the constant input and the noise) (hence its value at the next time step is different from zero in general).
 If a neuron does not fire and
does not receive influences from other neurons or input, then its membrane potential decays exponentially
fast with a decay rate $0 < \gamma< 1$.
The discussion of the biological relevance of this model and its extensions towards more
elaborated models with adaptive conductances has been done in \cite{cessac-vieville:08}.

\ssu{Technical definitions.}

\sssu{Spiking sequences.}
Call $\cM=\setR^N$ the phase space of our dynamical system. 
Given two integers $s<t$ (possibly negative) we note $\seq{V}{s}{t}$ the piece
of trajectory  $V(s), \dots, V(t)$.
To each membrane potential value, $V_i(t)$, we associate a variable $\omega_i(t)=Z(V_i(t))$.
The ``spiking pattern'' of the neural network at time $t$ is the vector 
$\omega(t)=\left(\omega_i(t)\right)_{i=1}^N$: it tells us which neurons are firing at time $t$,
($\omega_i(t)=1$) and which neurons are not firing at time $t$ ($\omega_i(t)=0$). 
We denote by $\seq{\omega}{s}{t}$ the sequence
or \textit{spike block } $\omega(s) \dots \omega(t)$. Associated with  each piece of trajectory  
$\seq{V}{s}{t}$ there is a unique spike block
  $\seq{\omega}{s}{t}$ with $\omega_i(n)=Z(V_i(n)), \ i=1\dots N, s \leq n \leq t$.
We note $Z(\seq{V}{s}{t})=\seq{\omega}{s}{t}$.
Also we note $\seq{\omega}{s}{t_1}\seq{\omega}{t_1}{t}=\seq{\omega}{s}{t}$ the concatenation
of the blocks $\seq{\omega}{s}{t_1}$ and $\seq{\omega}{t_1}{t}$.

\sssu{Raster plots.} 
Call $\cA$ the set of spiking patterns 
(alphabet). An element of  $\cA^\setZ$,
i.e. a bi-infinite sequence $\tom=\left\{\omega(t)\right\}_{t=-\infty}^{+\infty}$ 
of spiking patterns, is called a 
 ``raster plot''. It tells us which neurons are firing at each 
time $t \in \setZ$. In  experiments raster plots are
obviously finite sequences of spiking pattern but the extension
to $\setZ$, especially the possibility of considering an arbitrary
distant past (negative times) is quite useful in the present
work.
The set $\cA^\setZ$ is a topological space for the product topology \cite{kitchens:98}.
The open sets are the cylinder sets, namely the sets $[\seq{\omega}{s}{t}]=
\left\{\omega'  \in \cA^\setZ, \ \omega'(n)=\omega(n), \ n=s, \dots, t \right\}$.   
Cylinder sets are also a countable basis for the $\sigma$-algebra in $\cA^\setZ$. There is a natural distance on $\cA^\setZ$,
\beq\label{dX}
d_\Theta(\omega,\omega')=
\left\{
\baR{lll}
\left({\Theta^N}\right)^n, &\quad \mbox{if \, $\omega$ \, and \, $\omega'$ \, differ \, for
\, the \, first \, time \, in \, the \, $n$-th \, spiking \, pattern;}\\
0, &\quad \mbox{if \, $\omega=\omega'$},
\eaR
\right.
\eeq
for some $0 < \Theta < 1$. A classical choice
is $\Theta=\frac{1}{2}$. Here, it can be convenient to take $\Theta=\gamma$.

%

\sssu{Last firing time.} For $(s,t)\in \setZ^2,  s < t$, and each $i=1 \dots N$, we define the ``last firing
time of neuron $i$ in the sequence $\seq{\omega}{s}{t}$'' by: 
\beq\label{last_firing_time_s}
\tau_i(\seq{\omega}{s}{t})
\deq
\left\{
\baR{lll}
s, \quad \mbox{if} \quad \omega_i(k)=0, \quad k=s, \dots, t;\\
\max \left\{s \leq k \leq t,\omega_i(k)=1\right\}, \quad \mbox{if} \quad \exists k \in \left\{s, \dots, t\right\} \quad 
\mbox{such \, that} \quad \omega_i(k)=1.
\eaR
\right.
\eeq

\nid Therefore, $\tau_i(\seq{\omega}{s}{t})=s$ either if neuron $i$ fires
at time $s$ or if it does not fire during the whole time interval $[s,t]$.
In this way, the name ``last firing time'' is a little bit confusing, but this has no incidence
on the mathematical developments.

\ssu{The asymptotic probability distribution of membrane potentials and raster plots.}

\sssu{Conditional probability distribution of $V(t+1)$.}

Call $P=\cN(0,1)^{\otimes N\setZ}$,
 the joint distribution of the noise trajectories. 
Under $P$ the membrane potential $V$ is a stochastic process whose evolution is given eq. (\ref{DNN}). 
Fix a pair of integers $(s,t)$, $s < t$.
The probability distribution of $V(t+1)$ can be explicitly obtained with the following
remark. Since the cylinder sets $[\seq{\omega}{s}{t}]$ constitute a (countable) basis for the $\sigma$-algebra
in $\cA^\setZ$ and since to each piece of trajectory $\seq{V}{s}{t}$ is associated
a unique sequence  $\seq{\omega}{s}{t}$, we consider first the probability distribution
of $V(t+1)$ \textit{conditioned} by $Z(\seq{V}{s}{t})=\seq{\omega}{s}{t}$
and by the initial condition $V(s)$, assumed here to be \textit{bounded}.
 Then, the following holds as easily checked with a few algebra:

\bp\label{pVin}

For each $(s,t)\in \setZ^2,  s < t$, conditionally to  $Z(\seq{V}{s}{t})=\seq{\omega}{s}{t}$,
and given $V(s)$,
\beq\label{Vin}
V_i(t+1)=
\left\{
\baR{lll}
\gamma^{t+1-s} V_i(s)+C_i(\seq{\omega}{s}{t})+\sigma_B \xi_i(\seq{\omega}{s}{t}), 
\quad \mbox{if \, neuron \, $i$ \, didn't \, fire \, in \, the \, time \, interval \,} $[s,t]$;\\
C_i(\seq{\omega}{s}{t})+\sigma_B \xi_i(\seq{\omega}{s}{t}), 
\quad \mbox{otherwise}.\\
\eaR
\right.
\eeq
where 
\beq\label{Cist}
C_i(\seq{\omega}{s}{t})=
\sum_{j=1}^N W_{ij}x_{ij}(\seq{\omega}{s}{t})+\Iei \frac{1-\gamma^{t+1-\tau_i({\seq{\omega}{s}{t}})}}{1-\gamma},
\eeq
\beq\label{xijst}
x_{ij}(\seq{\omega}{s}{t})=\sum_{l=\tau_i(\seq{\omega}{s}{t})}^{t} \gamma^{t-l}\omega_j(l),
\eeq
\beq\label{xiist}
\xi_i(\seq{\omega}{s}{t})=\sum_{l=\tau_i(\seq{\omega}{s}{t})}^{t} \gamma^{t-l}B_i(l).
\eeq
\ep

\textbf{Remark.}
\bit
\item This equation expresses that the neuron loses its memory whenever it fires. This is due to the fact that
we reset the membrane potential, after firing. This consequently
simplifies the following analysis. For a discussion on dynamics of spiking neural phase models
when the condition is relaxed see \cite{kirst-geisel-etal:09}.
\item Clearly, the membrane potential is the sum of a ``deterministic'' part, $
\gamma^{t+1-s} V_i(s)+C_i(\seq{\omega}{s}{t})$, fixed by initial condition
at time $s$ and by the spike sequence $\seq{\omega}{s}{t}$,
and a stochastic part, $\sigma_B \xi_i(\seq{\omega}{s}{t})$, where the probability distribution of
the noise $\xi_i(\seq{\omega}{s}{t})$ is also fixed by the spike sequence $\seq{\omega}{s}{t}$.
More precisely, since the $B_i$'s are independent, Gaussian with mean zero and variance $1$,
the  $\xi_i(\seq{\omega}{s}{t})$'s, $i=1 \dots N$ are, under $P$, Gaussian, independent,  
with zero mean and variance $\frac{1-\gamma^{2(t+1-\tau_i(\seq{\omega}{s}{t}))}}{1-\gamma^2}.$
\eit

Denote by $\Exp{}$ the expectation under $P$.
It follows that:

\bp\label{CondLawV(t+1)}
For each $(s,t)\in \setZ^2,  s < t$, conditionally to  $Z(\seq{V}{s}{t})=\seq{\omega}{s}{t}$,
and given $V(s)$, 
$V(t+1)$ is Gaussian with mean:
%
$$\Exp{V_i(t+1)\, | \,\seq{\omega}{s}{t},V(s) }=
\left\{
\baR{lll}
\gamma^{t+1-s} V_i(s)+C_i(\seq{\omega}{s}{t}), 
\quad \mbox{if \, neuron \, $i$ \, didn't \, fire \, in \, the \, time \, interval \,} $[s,t]$;\\
C_i(\seq{\omega}{s}{t}), 
\quad \mbox{otherwise}.\\
\eaR
\right.$$
%
and covariance:
%
$$\CovB{V_i(t+1),V_j(t+1)\, | \,\seq{\omega}{s}{t},V(s)}= \sigma^2_i(\seq{\omega}{s}{t})\delta_{ij}.$$
%
%
with:
\beq\label{sigmai_cond_Gen}
\sigma^2_i(\seq{\omega}{s}{t})=\sigma_B^2\frac{1-\gamma^{2(t+1-\tau_i(\seq{\omega}{s}{t}))}}{1-\gamma^2}.
\eeq

Thus, the $V_i(t+1)$'s, $i=1 \dots N$,  are conditionally independent.
\ep

\textbf{Remark} We used a slight abuse of notation since we condition by $\seq{\omega}{s}{t}$
instead of $Z(\seq{V}{s}{t})=\seq{\omega}{s}{t}$.

\sssu{The probability that  some neuron $i$ does not fire within the time interval $[s,t]$}

It is given by:
$$ P\left(\bigcap_{n=s}^t \left\{ V_i(n) \, < \, \theta\right\}\right)=
\sum_{\seq{\omega}{s}{t} \in \cA^{t-s}} 
P\left(\bigcap_{n=s}^t \left\{ V_i(n) \, < \, \theta\right\} \, | \, \seq{\omega}{s}{t}\right)
P\left( \seq{\omega}{s}{t} \right)
$$
$$
=
\sum_{\seq{\omega}{s}{t} \in \cA^{t-s}} 
\prod_{n=s+1}^t
P\left(\left\{ V_i(n) \, < \, \theta\right\} \, | \, 
\bigcap_{l=s}^{n-1} \left\{ V_i(l) \, < \, \theta\right\} \cap \seq{\omega}{s}{t}\right)
P\left( \left\{ V_i(s) \, < \, \theta\right\} \, | \, \seq{\omega}{s}{t} \right)
P\left( \seq{\omega}{s}{t} \right)
$$
$$
=
\sum_{\seq{\omega}{s}{t} \in \cA^{t-s}} 
\prod_{n=s+1}^t P\left( \left\{ V_i(n) \, < \, \theta\right\} \, | \, 
\bigcap_{l=s}^{n-1} \left\{ V_i(l) \, < \, \theta\right\} \cap \seq{\omega}{s}{n-1}\right)
P\left( \left\{ V_i(s) \, < \, \theta\right\} \, | \, \seq{\omega}{s}{s} \right)
P\left( \seq{\omega}{s}{t} \right)
$$

From prop. (\ref{CondLawV(t+1)}), we have:
$$P\left( \left\{ V_i(n) \, < \, \theta\right\} \, | \, 
\bigcap_{l=s}^{n-1} \left\{ V_i(l) \, < \, \theta\right\} \cap \seq{\omega}{s}{n-1}\right)
=
P\left( \gamma^{n-s} V_i(s)+C_i(\seq{\omega}{s}{n-1})+\sigma_B \xi_i(\seq{\omega}{s}{n-1}) \, < \, \theta \, 
\right),$$
where $C_i(\seq{\omega}{s}{n-1})$ and $\xi_i(\seq{\omega}{s}{n-1})$ are given
by (\ref{Cist}),(\ref{xijst}),(\ref{xiist}) with $\tau_i(\seq{\omega}{s}{t})=s$.
Since, in this case, $\xi_i(\seq{\omega}{s}{n-1})$ is Gaussian, centered, with variance 
$ \frac{1-\gamma^{2(n-s)}}{1-\gamma^2}$ (eq. \ref{sigmai_cond_Gen}) we have:
$$P\left( \left\{ V_i(n) \, < \, \theta\right\} \, | \, 
\bigcap_{l=s}^{n-1} \left\{ V_i(l) \, < \, \theta\right\} \cap \seq{\omega}{s}{n-1}\right)
=1-\pi\left(\frac{\theta-\gamma^{n-s} V_i(s)-C_i(\seq{\omega}{s}{n-1})}{\sigma_B \sqrt{\frac{1-\gamma^{2(n-s)}}{1-\gamma^2}} } \right).$$

Since  $V_i(s)$ and
the $W_{ij}$'s are assumed to be bounded we have, whatever $n>s$:
\beq
0 <\Pi_- < P\left( \left\{ V_i(n) \, < \, \theta\right\} \, | \, 
\bigcap_{l=s}^{n-1} \left\{ V_i(l) \, < \, \theta\right\} \cap \seq{\omega}{s}{n-1}\right)
< \Pi_+ < 1, 
\eeq
for some constants $\Pi_-,\Pi_+$ depending on parameters $\gamma,W_{ij}, \, i,j =1 \dots N, I_i, \, i=1 \dots N$.
Likewise,  $0 < a < P\left( \left\{ V_i(s) \, < \, \theta\right\} \, | \, \seq{\omega}{s}{s} \right) < b <1$.
 Without loss of generality, e.g. redefining $\Pi_-$ as 
$\min(\Pi_-,a)$ (redefining $\Pi_+$ as $\max(\Pi_+,b)$) we may write 
$0 < \Pi_- < P\left( \left\{ V_i(s) \, < \, \theta\right\} \, | \, \seq{\omega}{s}{s} \right) < \Pi_+ <1$.
As a consequence, 

\bp\label{boundsPtau}
The  probability that  some neuron $i$ does not fire within the time interval $[s,t]$ has
the following bounds:
$$0 < \Pi_-^{t-s} <
P\left(\bigcap_{n=s}^t \left\{ V_i(n) \, < \, \theta\right\}\right)
< \Pi_+^{t-s} < 1.$$
As a consequence, whatever $s < t$, $t-s$ finite, there is a positive probability
that  some neuron $i$ does not fire within the time interval $[s,t]$.
This probability vanishes exponentially fast as $|t-s| \to +\infty$.
\ep

\sssu{Permanent regime}\label{S_perm_reg}

The main drawback of the previous results is that we have to condition on 
the ``initial'' condition $V(s)$ for spiking sequences such that some neuron does not fire
between $s$ and $t$. But the probability distribution of $V(s)$ is not
known. It has either to be ``guessed'' from ad hoc assumptions: is it Gaussian, uniform, ``fractal''  ... ?
Actually, the determination of initial conditions distribution is, to our opinion, one of the main
obstacle toward realistic characterizations or simulations of neural network models, intended
to somehow mimics the dynamics (of some part) of the brain, at some stage of its evolution.
Indeed, when considering the evolution of a set of neurons, one starts from
some ``initial'' time $s$ which corresponds to the beginning of the experiment.
This is NOT the beginning of the system under study, which has undergone a previous
evolution that actually determines the distribution of membrane potentials
at time $s$. This distribution has little chances to be Gaussian or anything so mathematically ``convenient'',
unless one finds strong arguments to justify this. Actually, as we show,
such assumption is wrong in model (\ref{DNN}). Therefore, to compute 
the  distribution of membrane potential
at time $s$ one has to consider the previous evolution of the system, which only postpones the problem,...
unless one assumes that this initial condition was drawn in an \textit{infinite past}.
On phenomenological grounds, ``infinite past'' means ``a time quite longer
than all characteristic time scales in the system'', though, mathematically,
one may take it truly infinite.
This is what we do here, focusing on what we call a ``permanent
regime'' (by analogy with Physics) where the initial condition is fixed in the infinite
past, namely $s \to -\infty$. In this case, indeed, $\gamma^{t+1-s} V_i(s) \to 0$.
As we show, this procedure selects a \textit{unique} probability distribution
for membrane potentials, with a highly non trivial structure
(see eq. \ref{rhoiV}).

We therefore consider left-infinite sequences $\seq{\omega}{-\infty}{t}$ with corresponding last firing time: 
\beq\label{last_firing_time}
\tau_i(\seq{\omega}{-\infty}{t})
\deq
\left\{
\baR{lll}
-\infty, \quad \mbox{if} \quad \omega_i(k)=0, \quad \forall k \leq t;\\
\max \left\{-\infty < k \leq t,\omega_i(k)=1\right\} \quad \mbox{otherwise}.
\eaR
\right.
\eeq

We now show that proposition \ref{CondLawV(t+1)} extends as well to the case $s \to -\infty$,
namely:

\bp\label{AsCondLawV(t+1)}
For each $t \in \setZ$, conditionally to  $\seq{\omega}{-\infty}{t}$,
$V(t+1)$ is Gaussian with mean:
\beq\label{Cit}
\Exp{V_i(t+1)\, | \,\seq{\omega}{-\infty}{t} }=C_i(\seq{\omega}{-\infty}{t})=
\sum_{j=1}^N W_{ij}x_{ij}(\seq{\omega}{-\infty}{t})+\Iei \frac{1-\gamma^{t+1-\tau_i({\seq{\omega}{-\infty}{t}})}}{1-\gamma},
\eeq
where,
\beq\label{xijt}
x_{ij}(\seq{\omega}{-\infty}{t})=\sum_{l=\tau_i(\seq{\omega}{-\infty}{t})}^{t} \gamma^{t-l}\omega_j(l),
\eeq
and covariance:
%
\beq
\CovB{V_i(t+1),V_j(t+1)\, | \,\seq{\omega}{- \infty}{t}}=\sigma^2_i(\seq{\omega}{-\infty}{t})\delta_{ij}=\sigma_B^2\frac{1-\gamma^{2(t+1-\tau_i(\seq{\omega}{-\infty}{t}))}}{1-\gamma^2}\delta_{ij}.
\eeq
Thus, the $V_i(t+1)$'s, $i=1 \dots N$,  are conditionally independent.
\ep
  
\bpr
Let us first show that the quantities defined by eq. (\ref{Cist}),(\ref{xijst}),(\ref{xiist}) are well
defined in the limit $s \to -\infty$. Consider first the limit of $x_{ij}(\seq{\omega}{s}{t})$
given by eq. (\ref{xijst}). There are two possibilities. Either $\seq{\omega}{-\infty}{t}$
is such that $\tau_i(\seq{\omega}{-\infty}{t})=n > -\infty$. Then,  $x_{ij}(\seq{\omega}{-\infty}{t})$
is a finite sum $\sum_{l=n}^{t} \gamma^{t-l}\omega_j(l)$ and is well defined.
Or, $\tau_i(\seq{\omega}{-\infty}{t})= -\infty$. Then 
$x_{ij}(\seq{\omega}{-\infty}{t})=\sum_{l=-\infty}^{t} \gamma^{t-l}\omega_j(l)=
\sum_{l=0}^{+\infty} \gamma^{l}\omega_j(t-l)$. This series  converges
since $\gamma<1$ and $\omega_j(t-l) = 0,1$. 
Moreover,
$\xi_i(\seq{\omega}{s}{t})=\sum_{l=\tau_i(\seq{\omega}{s}{t})}^{t} \gamma^{t-l}B_i(l)$
is a sum (possibly infinite) of independent Gaussian centered variables with finite variance.
As a consequence,
\beq\label{xiit}
\xi_i(\seq{\omega}{-\infty}{t})=\sum_{l=\tau_i(\seq{\omega}{-\infty}{t})}^{t} \gamma^{t-l}B_i(l),
\eeq
is Gaussian centered with variance $\frac{1-\gamma^{2(t+1-\tau_i(\seq{\omega}{-\infty}{t}))}}{1-\gamma^2}.$
Finally, from (\ref{Vin}) 
$V_i(t+1)=C_i(\seq{\omega}{-\infty}{t})+\sigma_B \xi_i(\seq{\omega}{-\infty}{t}),$
and the proposition follows from the independence of the $\xi_i(\seq{\omega}{-\infty}{t})$'s
and their Gaussian distribution.
\epr

\sssu{Elementary bounds.} 
We have:

\beq\label{boundxij_as}
0  \leq x_{ij}(\seq{\omega}{-\infty}{t}) \leq \frac{1}{1-\gamma},
\eeq
and, 
\beq\label{boundCi_as}
C_i^- \deq \Iei + \frac{1}{1-\gamma}\sum_{\stackrel{j=1}{W_{ij}<0}}^N W_{ij}
\leq C_i(\seq{\omega}{-\infty}{t}) 
\leq 
\frac{1}{1-\gamma}
\left(
\sum_{\stackrel{j=1}{W_{ij}>0}}^N W_{ij}+\Iei.
\right)\deq C_i^+,
\eeq
 In the same way 
\beq\label{boundsi_as}
\sB \leq \sigma^2_i(\seq{\omega}{-\infty}{t}) \leq  \frac{\sB}{1-\gamma^2}.
\eeq

\sssu{The transition probability.} 

We now compute the probability of a spiking pattern at time $t+1$,
$\omega(t+1)$, given the past $\seq{\omega}{-\infty}{t}$. It is given
by the following:

\bp
The probability  of $\omega(t+1)$ conditionally to  $\seq{\omega}{-\infty}{t}$
 is given by:
\beq\label{Pt+1_cond}
P\left(\omega(t+1)| \seq{\omega}{-\infty}{t},\right)=\prod_{i=1}^N P\left(\omega_i(t+1)| \seq{\omega}{-\infty}{t}\right),
\eeq

\nid with

\beq\label{Pt+1_cond_i}
 P\left(\omega_i(t+1)| \seq{\omega}{-\infty}{t},\right)=
\omega_{i}(t+1)\pi\left(\frac{\theta-C_i(\seq{\omega}{-\infty}{t})}{\sigma_i(\seq{\omega}{-\infty}{t})}\right)+
\left(1-\omega_{i}(t+1)\right)\left(1-\pi\left(\frac{\theta-C_i(\seq{\omega}{-\infty}{t})}{\sigma_i(\seq{\omega}{-\infty}{t})}\right)\right).
\eeq
\ep

\bpr
We have, using the conditional independence of the $V_i(t+1)$'s:
$$P\left(\omega(t+1) \, | \, \seq{\omega}{-\infty}{t}\right)=
%
%
\prod_{i=1}^N
\left[
\omega_i(t+1) \, 
P\left(
V_i(t+1) \geq \theta\, | \, \seq{\omega}{-\infty}{t} 
\right)
+
(1-\omega_i(t+1)) \, P\left(V_i(t+1)< \theta \, \Big| \, \seq{\omega}{-\infty}{t}
\right)
\right].
$$
Since the $V_i(t+1)$'s are Gaussian, with mean $C_i(\seq{\omega}{-\infty}{t})$
with a variance $\sigma^2_i(\seq{\omega}{-\infty}{t})$ we directly obtain (\ref{Pt+1_cond}),(\ref{Pt+1_cond_i}).
\epr

Consequently, it is possible, knowing the past sequence $ \seq{\omega}{-\infty}{t}$, to determine the probability
of the spiking pattern $\omega(t+1)$. In this way, $P\left(\omega(t+1)| \seq{\omega}{-\infty}{t},\right)$
acts as a transition probability, as in Markov chains. But here, the length of the Markov chain depends
on the last firing time of each neuron, since in fact,
$$P\left(\omega_i(t+1)| \seq{\omega}{-\infty}{t}\right)=
P\left(\omega_i(t+1)| \seq{\omega}{\tau_i(\seq{\omega}{-\infty}{t})}{t}\right).$$

The problem here is that, according to prop. \ref{boundsPtau},  we cannot \textit{bound} $\tau_i(\seq{\omega}{-\infty}{t})$.
Although 
this time $\tau_i(\seq{\omega}{-\infty}{t})$ 
is almost-surely finite, nevertheless, whatever $T>0$, there is
a positive probability set of sequences $\omega$ such that 
$\tau_i(\seq{\omega}{-\infty}{t})< t-T$. So we have to consider
a process where transition probability may have an \textit{unbounded memory}.
This type of process is called ``variable length Markov chain'' \cite{maillard:07}.
Such processes  can be studied in the general
context of \textit{chains with complete connections} and  \textit{$g$-measures},
developed in section \ref{SGibbs}.

\sssu{Stationarity.}\label{SStat}

In the present setting where $I_i$ does not depend on $t$ we have the following
property: 

\bp\label{PStat}
Fix a sequence
$\seq{a}{-\infty}{0}$, $a(-n) \in \cA, \, n \geq 0$. 
Then, $\forall t \in \setZ$,
\beq\label{EStat}
P\left(\omega(t)=a(0) \, | \, \seq{\omega}{-\infty}{t-1}=\seq{a}{-\infty}{-1} \right)=
P\left(\omega(0)=a(0) \, | \, \seq{\omega}{-\infty}{-1} =\seq{a}{-\infty}{-1}  \right).
\eeq
\ep

\bpr
Assume that $\omega(t-n)=a(-n), n \geq 0$, as in the l.h.s of (\ref{EStat}). 
Then, according to eq. (\ref{last_firing_time}),
$\tau_i(\seq{\omega}{-\infty}{t-1})=t+\tau_i(\seq{a}{-\infty}{-1})$. 
 Therefore,
according to eq. (\ref{xijt}),
$$x_{ij}(\seq{\omega}{-\infty}{t-1})=\sum_{l=\tau_i(\seq{\omega}{-\infty}{t-1})}^{t-1} \gamma^{t-1-l}\omega_j(l)
=\sum_{l=t+\tau_i(\seq{a}{-\infty}{-1})}^{t-1} \gamma^{t-1-l}a_{j}(l-t)
=\sum_{l'=\tau_i(\seq{a}{-\infty}{-1})}^{-1} \gamma^{-1-l'}a_{j}(l')=x_{ij}(\seq{a}{-\infty}{-1}),$$
and, according to (\ref{Cit}),
$$C_i(\seq{\omega}{-\infty}{t-1})=
\sum_{j=1}^N W_{ij}x_{ij}(\seq{\omega}{-\infty}{t-1})+\Iei \frac{1-\gamma^{t-\tau_i({\seq{\omega}{-\infty}{t-1}})}}{1-\gamma},
=\sum_{j=1}^N W_{ij}x_{ij}(\seq{a}{-\infty}{-1})+\Iei \frac{1-\gamma^{-\tau_i({\seq{a}{-\infty}{-1}})}}{1-\gamma}
=C_i(\seq{a}{-\infty}{-1}).
$$
Note that this last property holds because $I_i$ does not depend on time. We have also, from the same arguments,
$$
\sigma^2_i(\seq{\omega}{-\infty}{t-1})=
\sigma^2_i(\seq{a}{-\infty}{-1}).
$$
Consequently,
$$
P\left(\omega(t)=a(0) \, | \, \seq{\omega}{-\infty}{t-1}=\seq{a}{-\infty}{-1} \right)=
\prod_{i=1}^N 
\left[\omega_{i}(t) \, \pi\left(\frac{\theta-C_i(\seq{\omega}{-\infty}{t-1})}{\sigma_i(\seq{\omega}{-\infty}{t-1})}\right)+
\left(1-\omega_{i}(t)\right)\left(1-\pi\left(\frac{\theta-C_i(\seq{\omega}{-\infty}{t-1})}{\sigma_i(\seq{\omega}{-\infty}{t-1})}\right)\right)
\right]=
$$
$$
=\prod_{i=1}^N 
\left[a_{i}(0) \, \pi\left(\frac{\theta-C_i(\seq{a}{-\infty}{-1})}{\sigma_i(\seq{a}{-\infty}{-1})}\right)+
\left(1-a_{i}(0)\right)\left(1-\pi\left(\frac{\theta-C_i(\seq{a}{-\infty}{-1})}{\sigma_i(\seq{a}{-\infty}{-1})}\right)
\right)
\right]
=P\left(\omega(0)=a(0) \, | \, \seq{\omega}{-\infty}{-1} =\seq{a}{-\infty}{-1}  \right).
$$
\epr
 
Therefore, instead of considering a family of
transition probabilities depending on $t$, it suffices to
define the transition probability 
at one time $t \in \setZ$, for example  $t=0$.

\su{The equilibrium state.}\label{SGibbs}

In this section we show the existence of a unique invariant probability distribution for the dynamics (\ref{DNN}) and characterize it explicitly. Especially, we show that it is an \textit{equilibrium state} and a
 \textit{Gibbs state} in the sense of the thermodynamic formalism in ergodic theory \cite{keller:98,chazottes-keller:09}. 
For this we use the concept of $g$-measures \cite{keane:72}, coming from ergodic theory, and
 very close (equivalent
in the present setting) to the concept of chains with complete
connections\footnote{The concept of chains with complete connections, dates
 back to  1935 \cite{onicescu-mihoc:35}. They are a generalization
of Markov chains with an infinite memory. More precisely, they are induced by conditional probabilities of the form:
$P(\omega(t+1)  \, | \, \seq{\omega}{-\infty}{t}).$
These transition probabilities appear to be a extension of the notion of $k$-step Markov
chain with an infinite $k$. These objects must be taken with some precautions because, in
the non-Markovian case, the conditioning is always on an event of probability zero (see  \cite{maillard:07}
for proper definition).
The transition probabilities given by eq. (\ref{Pt+1_cond}),(\ref{Pt+1_cond_i}) define a system of such transition probabilities.}, 
which comes from probability theory. However, all known theorems used here being
formulated in the $g$-measure context, we use this formulation here.
For the convenience of the reader we however provide examples
and illustrations of the used notions.  Our main reference are \cite{chazottes:99,ledrappier:74,coelho-quas:98,bressaud-fernandez-etal:99,maillard:07}.

We proceed in several steps leading us to the  theorems  \ref{ThUnicity},\ref{eqstate}  the main results of this paper.

\ssu{Definitions and elementary results}

\sssu{Setting and notations.}

 The main object under study here is the family of transitions
probabilities (\ref{Pt+1_cond}). Using the stationarity proposition \ref{PStat} we can restrict to 
transition probabilities of the form $P(\omega(0)  \, | \, \seq{\omega}{-\infty}{-1})$.
Namely, we may focus on 
sequences in $\cA_{-\infty}^{0}$.
From now on we set 
$\omega=\seq{\omega}{-\infty}{0}$, $\uom=\seq{\omega}{-\infty}{-1}$ (a sequence $\uom$
is called an ``history''), 
$X=\seq{\cA}{-\infty}{0}$, $\uX=\seq{\cA}{-\infty}{-1}$. 
The $\sigma$-algebra (the set of cylinders) on $X$ (resp. $\uX$)
is denoted $\cF$ (resp. $\ucF$).

Call $T$ the \textit{right} shift over $X$ i.e. $(T\omega)(t)=\omega(t-1), \, t \leq 0$.
The use of the right shift, instead of the left shift currently used
in dynamical systems theory, is related to the formulation
of the problem in terms of transitions probabilities (see for example \cite{bressaud-fernandez-etal:99}).
We note $\omega a$, the right
concatenation of $\omega$ and $a \in \cA$, namely, this is the sequence 
$\omega'$ such that $\omega'(t-1)=\omega(t), \, t \leq 0 $ and $\omega'(0)=a$.
Note that $T(\omega a) = \omega$.

\sssu{$g$-functions.}

\bdf 
A \textit{$g$-function} over $(X,T)$ is a  measurable function $g : X  \to [0,1] $ which satisfies,
for all $\omega \in X$:
\beq\label{def_g_function}
\sum_{\omega', \, T(\omega')=\omega} g(\omega')=1.
\eeq
\edf

Examples of $g$ functions are precisely transition probabilities of type (\ref{Pt+1_cond}). 
Indeed, for $\omega \in X$, set: 
\beq\label{g0}
g_0(\omega)=
P\left(\omega(0) \, | \, \uom \right)
=\prod_{i=1}^N 
\left[
\omega_{i}(0)\pi\left(\frac{\theta-C_i(\uom)}{\sigma_i(\uom)}\right)+
\left(1-\omega_{i}(0)\right)
\left(
1-\pi\left(\frac{\theta-C_i(\uom)}{\sigma_i(\uom)}\right)
\right)
\right].
\eeq

Then, by definition, all $\omega'$'s in the sum (\ref{def_g_function}) have the form $\omega'=\omega a$,
and:
$$
\sum_{\omega', \, T(\omega')=\omega} g_0(\omega')=
\sum_{a \in \cA} 
 P\left(a \, | \, \uom \right)=1.
$$

We now give two properties of $g_0$ used below.

\sssu{$g_0$ is non-null}

A $g$ function is \textit{non null} on $X$ if
 for all $\omega \in X$, $g(\omega)>0$.
We have:

\bp\label{Pnonnull}
The $g$-function $g_0$, given by (\ref{g0}), is non-null.
\ep

\bpr
It suffices to check that  $P\left(\omega(0) \, | \, \uom \right) >0$.
 If  there exists
 $\omega \in X$ such that $ P\left(\omega(0) \, | \, \uom \right) = 0$,
then, for some $i \in \left\{1 \dots N \right\}$,
$\pi\left(\frac{\theta-C_i( \uom)}{\sigma_i( \uom)}\right)=0$
or $1$. 
This imposes that either $C_i( \uom)=\pm \infty$ or $\sigma_i( \uom)=0$
which is not possible since these quantities are bounded by bounds (\ref{boundCi_as}),(\ref{boundsi_as}).
\epr

\sssu{$g_0$ is continuous.}

\bdf The \textit{variation} of a  $g$-function $g$ is:
$$var_k(g)=\sup\left\{ \, | \,g(\omega)-g(\omega')\, | \,: \omega, \omega' \in X, \omega(t)=\omega'(t), \, \forall t \in \left\{-k, \dots, 0 \right\}\right\}.$$ 
\edf

\bdf
A $g$-function is \textit{continuous} if $var_k(g) \to 0$ as $k \to +\infty.$
\edf

\bp \label{go_cont}
$g_0$ is continuous.
\ep

\bpr 
We shall use the following inequalities. 

\ben
\item For a collection $0 \leq a_i,b_i \leq 1, \, \forall i=1 \dots N$,
we have\footnote{We thank one reviewer for this useful remark.} 
\beq\label{ineq1}
|\prod_{i=1}^N a_i-\prod_{i=1}^N b_i| \leq \sum_{i=1}^N |a_i- b_i|,
\eeq
as easily proved by recursion.

\item For $0 \leq x <1$, write $\sqrt{1-x} = 1-\sum_{n=1}^{+\infty} f_n x^n$, where 
$f_n={\frac {{4}^{-n} \left( 2\,n \right) !}
{ \left( n! \right) ^{2} \left( -1+2\,n \right) }} \geq 0
$
are the series coefficients of $\sqrt{1-x}$.
Then, for $A,B$ real, $0 \leq u,v < 1$,

\beq\label{ineq2}
|A\sqrt{1-u}-B\sqrt{1-v}| \leq | A - B| + \sum_{n=1}^{+\infty} f_n |A u^n - B v^n |
\leq | A - B| + \sum_{n=1}^{+\infty} f_n (|A| u^n + |B| v^n ). 
\eeq
\een

Fix $i \in \left\{1, \dots, N \right\}$. Set $y_i=\frac{\theta-C_i(\uom)}{\sigma_i(\uom)}$, $y'_i=\frac{\theta-C_i(\uomp)}{\sigma_i(\uomp)}$,
$C_i=C_i(\uom),C'_i=C_i(\uomp),\sigma_i=\sigma_i(\uom),\sigma'_i=\sigma_i(\uomp)$, $\tau_i=\tau_i(\uom)$, $\tau'_i=\tau_i(\uomp)$ to alleviate notations in the proof.
We have, for $k>0$, 
$$var_k(g_0)=\sup\left\{ \, 
|
\prod_{i=1}^N a_i-
\prod_{i=1}^N b_i
 | \,: \omega, \omega' \in X, \omega(t)=\omega'(t), \, \forall t \in \left\{-k, \dots, 0 \right\}\right\}.$$
\nid where $a_i=\omega_{i}(0)\pi\left(y_i\right)+\left(1-\omega_{i}(0)\right)\left(1-\pi\left(y_i\right)\right)$,
$b_i=\omega_{i}(0)\pi\left(y'_i\right)+\left(1-\omega_{i}(0)\right)\left(1-\pi\left(y'_i\right)\right)$.
Moreover, since either $\omega_i(0)=0$ or $\omega_i(0)=1$, 
$|a_i- b_i| = |\pi(y_i)-\pi(y'_i)|$. Therefore, using inequality (\ref{ineq1}),

$$var_k g_0 \leq 
\sum_{i=1}^N \sup\left\{ |
\pi\left(y_i\right) - \pi\left(y'_i)\right)
\ |: \omega, \omega' \in X, \omega(t)=\omega'(t), 
\forall t \in \left\{-k, \dots 0 \right\}\right\}.$$

Fix $k > 0$. Then, for all $\omega$ such that $\tiz \in \left\{-k \dots 0\right\}$,
$\pi\left(y_i\right)= \pi\left(y'_i\right)$.
Therefore, the $\sup$ is realized for those $\omega,\omega'$ such that $\tau_i,\tau'_i < -k$. 
Since $\gamma<1$ we have therefore $\gamma^{-\tau} < \gamma^k$ and
$\frac{1}{1-\gamma^{-\tau}}<\frac{1}{1-\gamma^k}$, for $\tau=\tau_i,\tau'_i$.
We have $|\pi(y_i)-\pi(y'_i)| \leq |y_i-y'_i| \|\pi'\|_{\infty}$ with $\|\pi'\|_{\infty}=\frac{1}{\sqrt{2\pi}}$. 
Moreover, 
$|y_i-y'_i| \leq \theta \left|\frac{1}{\sigma_i}- \frac{1}{\sigma'_i}\right|+\left|\frac{C'_i\sigma_i-C_i\sigma'_i}{\sigma_i \sigma'_i} \right|$.

We have, 
$$\left|\frac{1}{\sigma_i}- \frac{1}{\sigma'_i}\right|=
\frac{\sqrt{1-\gamma^2}}{\sigma_B} \left|\frac{1}{\sqrt{1-\gamma^{-2\tau_i}}} \, - \, \frac{1}{\sqrt{1-\gamma^{-2\tau'_i}}}  \right|$$
$$
= \frac{\sqrt{1-\gamma^2}}{\sigma_B
\sqrt{1-\gamma^{-2\tau_i}}\sqrt{1-\gamma^{-2\tau'_i}}} 
\left|\sqrt{1-\gamma^{-2\tau'_i}} - \sqrt{1-\gamma^{-2\tau_i}}\right|$$
$$
\leq \frac{2 \sqrt{1-\gamma^2}}{\sigma_B (1-\gamma^{2k})} 
\gamma^{2k} S(\gamma),$$
with 
\beq\label{Sgamma}
S(\gamma)=\sum_{n=1}^{+\infty} f_n \gamma^{2k(n-1)},
\eeq
the last inequality coming
from (\ref{ineq2}).

Likewise,
$$\left|\frac{C'_i\sigma_i-C_i\sigma'_i}{\sigma_i \sigma'_i} \right|=
\frac{\sqrt{1-\gamma^2}}{\sigma_B
\sqrt{1-\gamma^{-2\tau_i}}\sqrt{1-\gamma^{-2\tau'_i}}} 
\left|C'_i \, \sqrt{1-\gamma^{-2\tau_i}} - C_i \, \sqrt{1-\gamma^{-2\tau'_i}}\right|
$$
$$
\leq 
 \frac{\sqrt{1-\gamma^2}}{\sigma_B (1-\gamma^{2k})}
\left(
\left|C'_i -C_i \right| + \sum_{n=1}^{+\infty} f_n \gamma^{2kn}(|C'_i|+|C_i|)
\right)
$$
$$
\leq 
 \frac{ \sqrt{1-\gamma^2}}{\sigma_B (1-\gamma^{2k})}
\left(
\left|C'_i -C_i \right| + 2\gamma^{2k} |C_i^+| S(\gamma)
\right)
$$
where $C_i^+$ is given by (\ref{boundCi_as}).

We have,
$$
\left|C'_i - C_i \right|=
\left|
\sum_{j=1}^N W_{ij}
\left(
\sum_{l=\tau_i}^{-1} \gamma^{-l-1}\omega_j(l)
-
\sum_{l=\tau'_i}^{-1} \gamma^{-l-1}\omega'_j(l)
\right)
+ 
\frac{\Iei}{1-\gamma}
(\gamma^{-\tau_i}-\gamma^{-\tau'_i})
\right|
$$
$$
\leq
\sum_{j=1}^N \left|W_{ij} \right|
\left|
\sum_{l=\tau_i}^{-1} \gamma^{-l-1}\omega_j(l)
-
\sum_{l=\tau'_i}^{-1} \gamma^{-l-1}\omega'_j(l)
\right|
+ 
\frac{\left|\Iei\right|}{1-\gamma}
\left|\gamma^{-\tau_i}-\gamma^{-\tau'_i}\right|.
$$

Remark that  
$$\sum_{l=-k}^{-1} \gamma^{-l-1}\omega_j(l)=
\sum_{l=-k}^{-1} \gamma^{-l-1}\omega'_j(l),$$
since $\omega_j(l)=\omega'_j(l), l=-k \dots -1$. Therefore,
$$
\left|
\sum_{l=\tau_i}^{-1} \gamma^{-l-1}\omega_j(l) 
-
\sum_{l=\tau'_i}^{-1} \gamma^{-l-1}\omega'_j(l)
\right|
=
\left|
\sum_{l=\tau_i}^{-k-1} \gamma^{-l-1}\omega_j(l) 
-
\sum_{l=\tau'_i}^{-k-1} \gamma^{-l-1}\omega'_j(l)
\right|
\leq
  \frac{2\gamma^k}{1-\gamma}.$$
Moreover,
$$
\frac{\left|\Iei\right|}{1-\gamma}
\left|\gamma^{-\tau_i}-\gamma^{-\tau'_i}\right|
\leq 
2\frac{\left|\Iei\right|}{1-\gamma}
\gamma^{k}
$$
Finally,
$$
\left|C'_i - C_i \right|
\leq
\frac{2\gamma^k}{1-\gamma}
\left(\sum_{j=1}^N |W_{ij}|
+\left|\Iei\right|
\right)
.$$

Summarizing,
$$
\left|\frac{C'_i\sigma_i-C_i\sigma'_i}{\sigma_i \sigma'_i} \right|
\leq
\gamma^k \frac{2 \sqrt{1-\gamma^2}}{\sigma_B (1-\gamma^{2k})}
\left(\frac{1}{1-\gamma}
\left(\sum_{j=1}^N |W_{ij}|
+\left|\Iei\right|\right)
 + \gamma^{k}|C_i^+| S(\gamma)
\right)
$$
$$\|\pi(y_i)-\pi(y'_i) \| \leq 
\sqrt{\frac{2}{\pi}}
 \frac{ \sqrt{1-\gamma^2}}{\sigma_B (1-\gamma^{2k})}
\left(\frac{1}{1-\gamma}
\left(\sum_{j=1}^N |W_{ij}|
+\left|\Iei\right|\right)
 +  \gamma^k(\theta +|C_i^+|) S(\gamma)
\right)
\gamma^k,
$$
and,
$$
var_k g_0 \leq 
\sqrt{\frac{2}{\pi}}
 \frac{ \sqrt{1-\gamma^2}}{\sigma_B (1-\gamma^{2k})}
\left(\frac{1}{1-\gamma}
\left(\sum_{i,j=1}^N |W_{ij}|
+\sum_{i=1}^N\left|\Iei\right|\right)
 +  \gamma^k (N \theta + \sum_{i=1}^N |C_i^+|) S(\gamma)
\right)\gamma^k,
$$

Remarking  that $S(\gamma) \to f_1$ as $k \to +\infty$ we conclude that
$var_k g_0 \to 0$ as $k \to +\infty$ and behaves like $K \gamma^k$,
where the constant 

\beq\label{K}
K= \sqrt{\frac{2}{\pi}} \frac{1}{\sigma_B} \sqrt{\frac{1+\gamma}{1-\gamma}}
\left[
\sum_{i,j=1}^N |W_{ij}|
+\sum_{i=1}^N\left|\Iei\right|
\right],
\eeq
 depends on model parameters (synaptic weights and current)
and $N$. 
\epr

\textbf{Remark.} Note that $g_0$ is therefore
H\"olderian for the metric (\ref{dX})
(with an exponent $\frac{\log\gamma}{\log\Theta}$ so it is Lipschitz
for $\Theta=\gamma$). 
 
\ssu{The Gibbs equilibrium state}

We now prove that the system (\ref{DNN}) admits a unique invariant
probability measure (also called a $g$-measure). This measure
satisfies a variational principle (equilibrium state) and
has the form of a Gibbs distribution in statistical physics. 

\sssu{$g$-measure.}

\bdf
Let $g$ be a $g$-function. A probability measure $\mu$ in $\cP(X,\cF)$ is a $g$-\textit{measure}
if :
$$\int f(\omega)g(\omega a) \mu(d\omega) = \int_{\left\{\omega(0)=a \right\}} f(\omega)\mu(d\omega),$$
 $\forall a \in \cA$ and $\forall f$ measurable with respect to $\ucF$.
\edf

\sssu{There exists a $g$-measure for (\ref{DNN}) and it is unique.}

Since $g_0$ is continuous there is always a $g$-measure.
Now, a theorem of Johansson and Oberg 
 \cite{johansson-oberg:2002} states that if
 $g_0$ is  a continuous,  non-null $g$ function on $X$ satisfying:
\beq
\sum_{k\geq 0} var^2_k(\log g_0)) < +\infty,
\eeq
then the $g$-measure is unique.

\bth\label{ThUnicity}
The dynamical system (\ref{DNN}) has a unique 
$g_0$-measure whatever the values of parameters
$W_{ij}, \, i,j=1 \dots N, I_i, \, i=1 \dots N, \gamma, \theta$.
\enth

\bpr
This follows from the theorem of Johansson and Oberg.
 Indeed, using a proof similar to
prop. \ref{go_cont} and the same notations, we have, for $k>0$,
and using that 
$$\log\left[\omega_{i}(0)\pi(y_i)+\left(1-\omega_{i}(0)\right)\left(1-\pi(y_i)\right)\right]=
\omega_{i}(0)\log\left(\pi(y_i)\right)+\left(1-\omega_{i}(0)\right)\log\left(1-\pi(y_i)\right),$$
$$
\baR{lll}
var_k(\log g_0) =\\
\sup\left\{ \, 
| \,
\sum_{i=1}^N 
\left[
\omega_{i}(0)\left[\log\left(\frac{\pi(y_i)}{\pi(y_i')}\right)\right]
+\left(1-\omega_{i}(0)\right) \left[\log\left(\frac{1-\pi(y_i)}{1-\pi(y'_i)}\right)\right]
\right]
\, | \,: \omega, \omega' \in X, \omega(t)=\omega'(t), \, \forall t \in \left\{-k, \dots, 0 \right\}
\right\}
\eaR
$$
$$\leq 
\sum_{i=1}^N 
\sup\left\{ |
\log\left(\frac{\pi(y_i)}{\pi(y_i')}\right) |
+ |
\log\left(\frac{1-\pi(y_i)}{1-\pi(y_i')}\right) 
\ |: \omega, \omega' \in X, \omega(t)=\omega'(t), 
\forall t \in \left\{-k, \dots 0 \right\}\right\}.$$
From the bounds (\ref{boundCi_as}), (\ref{boundsi_as}), 
$$\sqrt{1-\gamma^2} \, \frac{\theta - C_i^+}{\sigma_B} \leq y_i, y'_i \leq \frac{\theta - C_i^-}{\sigma_B}.$$
Set $a \deq \sqrt{1-\gamma^2} \, \min_{i=1 \dots N} \frac{\theta - C_i^+}{\sigma_B}$
and $b \deq \max_{i=1 \dots N} \frac{\theta - C_i^-}{\sigma_B} < \infty$.
Denote  $\| f \|_{_{[a,b]}} \deq \sup_{x \in [a,b]} \left|f(x)\right|$, we have:
$$
var_k(\log g_0) \leq 2\, \|\frac{\pi'}{\pi} \|_{_{[a,b]}}
 \sum_{i=1}^N 
\sup
\left\{ 
|y_i
-
 y'_i
\ |: \omega, \omega' \in X,
 \omega(t)=\omega'(t), 
\forall t \in \left\{-k, \dots 0 \right\}
\right\},
$$
where  the norm 
$$ \|\frac{\pi'}{\pi} \|_{_{[a,b]}} = \frac{e^{-\frac{b^2}{2}}}{\int_{b}^{+\infty} e^{-\frac{u^2}{2}}du},$$
is finite since $b$ is finite. For the term  
$\max_i 
\sup
\left\{ 
|y_i 
-
y'_i)
\ |: \omega, \omega' \in X,
 \omega(t)=\omega'(t), 
\forall t \in \left\{-k, \dots 0 \right\}
\right\}$, we have the same majoration as in the proof of 
prop \ref{go_cont}. Thus,
$$var_k(\log g_0) \leq  K' \gamma^k,$$
\nid with 
\beq\label{Kp}
K'=\sqrt{2\pi}\, \|\frac{\pi'}{\pi} \|_{_{[a,b]}} K,
\eeq
$K$ given by (\ref{K}).
It follows that $\sum_{k\geq 0} var^2_k(\log g_0) \leq K'^2 \sum_{k\geq 0}(\gamma^{2})^k < \infty. $
Then the $g_0$ measure is unique.
\epr

Let us now characterize the structure of this $g$-measure.

\sssu{The Ruelle-Perron-Frobenius operator.}

For the $g$-function $g_0$ define the \textit{transfer operator} or \textit{Ruelle-Perron-Frobenius operator} 
$\cL_{g_0}$ from $C(X,\setR)$ to $C(X,\setR)$,
where $C(X,\setR)$ is the set of continuous real functions on $X$, by:
\beq \label{RPF}
\cL_{g_0} f(\omega)=\sum_{\omega' \, : \, T(\omega')=\omega} g_0(\omega')f(\omega'). 
\eeq
Denoting
 $\cL_{g_0}^n$,
$n>0$, the $n$-iterates of the RPF operator, $\cL_{g_0}^n f$ is the conditional expectation of $f$ on the time interval 
$\left[0,n-1\right]$ given the history $\uom$.
The Ruelle-Perron-Frobenius is the extension of matrices of probability transitions
for Markov chains. 

The adjoint of $\cL_{g_0}$ maps the set of probability measures on $X$ to itself and is defined by:
%
$$\cL_{g_0}^\ast\mu(f)=\mu(\cL_gf).$$
%

A probability measure $\mu$ on $X$ is a $g$-measure if and only if
$\cL_{g_0}^\ast\mu=\mu$ \cite{ledrappier:74}.

\sssu{Equilibrium state.}

Let $\psi$ be a continuous function 
$X \to \setR$ such that $\sum_{k=0}^\infty var_k(\psi) <\infty$
(also called a \textit{regular potential} \cite{keller:98}). Call $P_T(X)$ the set of  $T$-invariant finite measures on $X$. 
For $\mu \in P_T(X)$ let
\beq\label{DefEntro}
h(\mu)=\limsup_{n \to + \infty} 
\frac{1}{n+1} \sum_{\cyl{\omega}{0}{n}} \mu(\cyl{\omega}{0}{n})\log\mu(\cyl{\omega}{0}{n}),
\eeq
be the entropy of $\mu$, where the sum holds over
all cylinders $\cyl{\omega}{0}{n}$ of length $n+1$. Note that the entropy can be defined
in a more general setting (see \cite{keller:98}). Here we take a definition
which corresponds more to the one used in neural networks dynamics
analysis. 

\bdf
An \textit{equilibrium state}, $\mu_\psi$, is a $T$-invariant 
 measure on $X$, such
that:
\beq\label{VarPrinc}
P(\psi) \deq h(\mu_\psi)+\mu_\psi(\psi)=\sup_{\mu \in P_T(X) } h(\mu)+\mu(\psi).
\eeq
\edf.

The quantity $P(\psi)$ is called the ``topological pressure''
\cite{sinai:72,bowen:98,ruelle:78}. This is a fundamental
quantity and we come back to it in section
\ref{SFiniteRange}. It is zero whenever the potential $\psi$ is normalized, which
is the case here since $\psi$ is the log of a conditional probability.

Ledrappier has shown \cite{ledrappier:74} that 
if $\psi$ is a regular potential then the equilibrium 
states for $\psi$ are the $g$-measures for a continuous $g$ function. 
In our case where the $g$-measure is unique,
 $g_0$ is related to the potential $\psi$ 
given by:
%
$$\psi(\omega) = \log g_0(\omega).$$
%

Therefore, the following holds.

\bth\label{eqstate}
Whatever the parameters values 
the system (\ref{DNN}) has a unique $g_0$ probability measure, $\mu_\psi$,
 which is an equilibrium state for the potential 
\beq\label{psi}
\psi(\omega) \equiv \psi(\seq{\omega}{-\infty}{0})
 =\log g_0(\omega)=\sum_{i=1}^N 
\left[\omega_{i}(0) 
\log\left(\pi\left(\frac{\theta-C_i(\uom)}
{\sigma_i(\uom)}\right)\right)+
\left(1-\omega_{i}(0)\right)
\log
\left(1-\pi\left(\frac{\theta-C_i(\uom)}{\sigma_i(\uom)}\right)\right)
\right].
\eeq
\enth

\ssu{Consequences.}

\sssu{Asymptotic distribution of membrane potentials}

A first consequence of proposition \ref{AsCondLawV(t+1)} and theorem \ref{ThUnicity}
is:

\bp
The membrane potential vector $V$ is stationary with a
product density $\rho_{V}(v)=\prod_{i=1}^N \rho_{_{V_i}}(v_i)$ where:
\beq\label{rhoiV}
\rho_{_{V_i}}(v)=
\int_{\uX}  
\frac{1}{\sqrt{2 \pi \sigma_i(\uom)}} \exp\left(-\frac{1}{2} \left(\frac{v-C_i(\uom)}{\sigma_i(\uom)} \right)^2 \right)
d\mu_\psi(\uom).
\eeq

Its expectation is $\mu_\psi\left[C_i(\uom) \right]$ and its variance $\mu_\psi\left[\sigma_i^2(\uom) \right]$. 
\ep
 
\textbf{Comments.} This density is a ``mixture'' of Gaussian densities, but it is not Gaussian. 
Each Gaussian density in the decomposition depends on a specific history $\uom$,
and the integral holds on the set of all possible histories with a
weight $\mu_\psi(\uom)$. Therefore, to obtain a closed form for the stationary density of membrane
potential we need to know the invariant probability $\mu_\psi$ which weights
the possible histories \textit{over an unbounded past}. It has therefore
a highly non trivial structure as announced in section \ref{S_perm_reg}.

\sssu{Firing rates.}

Define:
\beq \label{defcondrate}
r_i(\omega) \deq P(\omega_i(0)=1 \, | \, \uom)=
\pi\left[\frac{\theta-C_i(\uom)}{\sigma_i(\uom)} \right],
\eeq
the probability that neuron $i$ fires at time $0$
given the past $\uom$ and,
\beq \label{defrate}
r_i \deq \mu_\psi(\omega_i(0)=1),
\eeq
called ``the firing rate'' of neuron $i$.
We have:
\beq\label{Explicit_ri}
r_i=\mu_\psi
\left(r_i(\omega)\right)=
\mu_\psi
\left(
\pi\left[\frac{\theta-C_i(\uom)}{\sigma_i(\uom)} \right]
\right).
\eeq

\sssu{Entropy.}\label{Sentro}

It results from (\ref{VarPrinc}) and the normalization
of the potential $\psi$ that $0=h(\mu_\psi)+\mu_\psi(\psi)$.
Therefore:
$$
h(\mu_\psi)=
-\sum_{i=1}^N 
\mu_\psi\left(\omega_i(0)\,\log(\pi\left[\frac{\theta-C_i(\uom)}{\sigma_i(\uom)}\right]) \,
+
 \, (1-\omega_i(0)) \, \log(1-\pi\left[\frac{\theta-C_i(\uom)}{\sigma_i(\uom)}\right])\right).
$$
Since either $\omega_i(0)=0$ or $1$, we have:
$$
h(\mu_\psi)=-\sum_{i=1}^N 
\left[
r_i \,
\mu_\psi\left(\log(\pi\left[\frac{\theta-C_i(\uom)}{\sigma_i(\uom)}\right])\right) 
+
(1-r_i) \, 
\mu_\psi\left(\log(1-\pi\left[\frac{\theta-C_i(\uom)}{\sigma_i(\uom)}\right])\right)
\right],
$$
and finally,
\beq\label{h_explicit}
h(\mu_\psi)=-\sum_{i=1}^N 
\left[
r_i \,
\mu_\psi\left(\log r_i(\omega)\right) 
+
(1-r_i) \, 
\mu_\psi\left(\log(1-r_i(\omega))\right)
\right].
\eeq
This looks like the classical entropy for a Bernoulli
scheme but with a crucial difference: one has to take
the expectation of the $\log$ of the probability
instead of the $\log$ of the expectation.

Moreover $\psi(\omega) <0$ (the strict inequality comes from proposition 
\ref{Pnonnull}). Therefore, $\mu_\psi(\psi) < 0$ and:

\bp
The entropy $h(\mu_\psi)$, given by (\ref{h_explicit}), is  positive whatever the value of parameters
$W_{ij}, i,j=1 \dots N, I_i, i=1 \dots N,\theta,\gamma$.
\ep

Though this result appears ``evident'' a priori, it is proved here, as an easy consequence of the
variational principle (\ref{VarPrinc}). Moreover, it provides an explicit value for the entropy,
which depends on parameters.

\sssu{Gibbs state.}

In the present setting equilibrium states are Gibbs states \cite{keller:98}.
A Gibbs state for the potential
$\psi$ is a probability measure 
$\mu_\psi$ such
that  one can find some constants 
$P(\psi),c_1,c_2$ with $0 < c_1 \leq 1 \leq c_2$ such
that for all $n \geq 0$ and for all $\omega \in X$:
\beq\label{dGibbs}
c_1 \leq \frac{\mu_\psi\left(\cyl{\omega}{0}{n}\right)}
{\exp \left[ -(n+1) P(\psi)+\sum_{k=0}^n\psi(T^k\omega) \right]}  \leq c_2.
\eeq

Basically, the condition (\ref{dGibbs}) expresses that the measure
of the cylinder $\cyl{\omega}{0}{n}$ behaves like: 
\beq\label{Gibbs_def}
\mu_\psi\left(\cyl{\omega}{0}{n}\right) 
\sim \frac{
\exp{\sum_{k=0}^n\psi(\seq{\omega}{-\infty}{k})}
}
{Z^{(n+1)}_\psi(\omega)},
\eeq
\nid which has therefore the classical form of Gibbs distribution where spins chains
are replaced by sequences of spiking patterns and
where the normalization factor $Z^{(n+1)}_\psi(\omega)$ is analog to a ``partition
function'' (but depends on $\omega$). Note that $P(\psi)=\lim_{n \to +\infty} \frac{1}{n+1}\log Z^{(n+1)}_\psi(\omega)$,
where the limit exists and is constant for $\mu_\psi$-almost-every
$\omega$. Thus the topological pressure $P(\psi)$ is analog to a free
energy density.

\sssu{Kullback-Leibler divergence.}

Let $\mu,\nu$ be two $T$-invariant measures. 
The Kullback-Leibler divergence between $\mu$
and $\nu$ is given by:

\beq\label{dKL}
d(\mu,\nu)=
\limsup_{n \to +\infty} 
\frac{1}{n+1}\sum_{\cyl{\omega}{0}{n}} 
\mu\left(\cyl{\omega}{0}{n}\right)
\log\left[\frac{\mu\left(\cyl{\omega}{0}{n}\right)}
{\nu\left(\cyl{\omega}{0}{n}\right)} \right],
\eeq
where the sum holds on all possible cylinders
$\cyl{\omega}{0}{n}$. It
provides some notion of asymmetric ``distance'' between $\mu$ and $\nu$.
Minimizing this divergence,
corresponds to minimizing 
``what is not explained in the empirical measure $\mu$
by the theoretical measure $\nu$''.

The following holds. For $\mu$ an ergodic measure 
and $\mu_\psi$ a Gibbs state with a potential
$\psi$, both defined on the same set of sequences, one has 
\cite{bowen:98,ruelle:69,keller:98,chazottes-keller:09}:  
\beq\label{dKLpres}
d\left(\mu,\mu_\psi \right) = P(\psi) - \mu(\psi) - h(\mu).
\eeq

This result is used in the next section. 

\su{Finite range approximations.}\label{SFiniteRange}

\ssu{Constructing a Markov chain with finite memory.} \label{SMarkov}

The main difficulty in handling the transition probabilities (\ref{Pt+1_cond}) and the related equilibrium state
is that they depend on an history dating back to 
 $\tau_i(\seq{\omega}{-\infty}{t})$, where $\tau_i(\seq{\omega}{-\infty}{t})$ is unbounded.
On the other hand, the influence of 
the activity of the network, say  at time $-l$, on the 
membrane potential $V_i$ at time $0$, appearing
in the term 
$x_{ij}(\seq{\omega}{-\infty}{0})=\sum_{l=\tau_i(\seq{\omega}{-\infty}{0})}^{0} \gamma^{-l}\omega_j(l),$
(eq. \ref{xijt}) decays exponentially fast as $l \to -\infty$.
 Thus, one
may argue that after a characteristic time depending
on $\frac{1}{|\log(\gamma)|}$
the past network  activity has little influence on $V_i(0)$.
 We now make
this statement precise, especially evaluating the error attached to this
approximation, before exploring its consequences in section \ref{SSpikeStat}.

\sssu{Range-$R$ approximation.}

Assume that we want to approximate the statistics of
spikes, given by the dynamics (\ref{DNN}), by  fixing a finite time horizon
$R$ such that  the membrane potential at time $0$  depends on the past only up to some
finite time $-R$. 
In this way, we truncate the histories and we approximate the transition probabilities 
$P\left(\omega(0)\, | \, \seq{\omega}{-\infty}{-1}\right)$,
 with unbounded memory,
by transition probabilities $P\left(\omega(0)\, | \, \seq{\omega}{-R}{-1}\right)$
where $\tau_i(\seq{\omega}{-\infty}{-1})$ is replaced by $\tau_i(\seq{\omega}{-R}{-1})$
(see eq. (\ref{last_firing_time_s})), thus limiting memory to at most $R$ time steps in the past.
%
%
%
These approximated transition probabilities
 constitute therefore a Markov chain with a memory depth $R$. 
How good is this approximation ?
To answer this question let us first construct
the Markov chain within more details.

\sssu{Blocks coding.} 

Since we are now only 
considering finite histories given by 
spike blocks of length $R$, of the form $\seq{\omega}{-R}{-1}$,
we may encode each of these blocks  by an integer
\beq\label{w_coding}
w=\sum_{i=1}^N\sum_{n=-R}^{-1} 2^{(i-1)+(n+R)N}\omega_i(n).
\eeq
We write $w \sim \seq{\omega}{-R}{-1}$. 
These integers
or  \textit{words} constitute the states of the Markov chain.
We note 
$$\Omega^R \deq  \left\{0, \dots, 2^{NR}-1 \right\},$$
the set of words.

\sssu{Transition matrix.} \label{STrans}

From prop. \ref{PStat}, this chain is homogeneous,  i.e. transition probabilities does not
depend on time. They are  encoded in a
$2^{NR} \times 2^{NR}$  matrix $\cL^{(R)}$ with entries:
\beq\label{cL}
\cL^{(R)}_{w',w}\deq 
\left\{
\baR{ccc}
P\left(\omega(0)\, | \, \seq{\omega}{-R}{-1}\right), \, &\mbox{if}& 
w' \sim \seq{\omega}{-R}{-1}, w \sim \seq{\omega}{-R+1}{0},\\
0, \, &\mbox{otherwise}.& 
\eaR
\right.
\eeq

If $w' \sim \seq{\omega}{-R}{-1}, w \sim \seq{\omega}{-R+1}{0}$ we say that $w$ ``follows''
$w'$ or that the transition $w' \to w$ is ``legal''. 
Note that when using a matrix representation where
$w' \in \Omega^R$, $w \in \Omega^R$ not all transition
are legal ($w',w $ must correspond to overlapping blocks).
We therefore use  the convention that
 non-legal transitions have a zero probability. 

It results from prop. \ref{Pnonnull} that 
any legal transition has a positive probability.
Note that  the transition
matrix $\cL^{(R)}$ corresponds to the matrix representation of the Ruelle-Perron-Frobenius (\ref{RPF})
 operator
in the case of a finite memory.  

\sssu{Incidence matrix.}

 The set of transitions is encoded in a
$2^{NR} \times 2^{NR}$ matrix $\cI$, called \textit{incidence matrix}
 with entries:
\beq
\cI_{w'w}=
\left\{
\baR{ccc}
1, & \mbox{if \, $w' \to w$ \, is \, legal}&\\
0; & \mbox{otherwise.}&
\eaR
\right.
\eeq   
It is easy to show that $\cI$ is \textit{primitive} namely  
$\exists m>0$ such that $\forall w,w' \in {\Omega^R} \times {\Omega^R}$, $\cI^m_{w'w}>0$, where
$\cI^m$ is the $m$-th power of $\cI$. Indeed, having
$\cI^m_{w'w}>0$ means that there exists a raster plot
$\omega$ which contains the block $w'$ and the block
$w$ where the first spiking pattern of each block 
is separated by $m$ time steps. Therefore,
taking $m = R+1$ any raster containing the concatenation
of blocks $w'w$ satisfies the requirement.

\sssu{Range $R+1$ potential.} Using the same
representation as (\ref{w_coding}), but with
blocks of size $R+1$,  to each block $\seq{\omega}{-R}{0}$ of length $R+1$ 
we associate a word $W=\sum_{i=1}^N\sum_{n=-R}^{0} 2^{(i-1)+(n+R)N}\omega_i(n)
 \sim \seq{\omega}{-R}{0}$ and define:
\beq\label{psi_R}
\psi^{(R)}(W) =
\sum_{i=1}^N 
\left[\omega_{i}(0) 
\log
\left(
\pi
\left(
\frac{\theta-C_i(\seq{\omega}{-R}{-1})}{\sigma_i(\seq{\omega}{-R}{-1})}
\right)
\right)+
\left(1-\omega_{i}(0)\right)
\log\left(1-\pi\left(\frac{\theta-C_i(\seq{\omega}{-R}{-1})}{\sigma_i(\seq{\omega}{-R}{-1})}\right)\right)
\right],
\eeq
 called a \textit{range-$R+1$ potential}. It corresponds to an approximation
of the potential (\ref{psi}) when the memory depth of the chain is $R$. Then:
%
$$\cL^{(R)}_{w',w}=e^{\psi^{(R)}(W)}\cI_{w'w}.$$

\sssu{The Perron-Frobenius theorem.}\label{SPR}

Since $\cI$ is primitive and since all legal transitions
have a positive probability,  $\cL^{(R)}$ is primitive
and the Perron-Frobenius theorem holds
\cite{gantmacher:66,seneta:06}. $\cL^{(R)}$ has
a real positive eigenvalue $s$ with maximal modulus, 
isolated from the rest of the spectrum. Moreover, since $\cL^{(R)}$
is a transition probability, $s=1$.   The quantity 
$$P(\psi^{(R)})=\log s=0$$
is the topological pressure of the potential $\psi^{(R)}$ \cite{keller:98}.

The corresponding left and right eigenvectors are respectively denoted $l$ and $r$ i.e. 
$l \, \cL^{(R)} = s  \, l$ and $\cL^{(R)}  \, r = s   \, r$ where $r$ as positive
entries $r_w >0$. Without of generality we may assume that
$\langle l,r \rangle=1$, where $\langle \, , \rangle$ denotes
the standard scalar product. 

The Markov chain has a unique invariant probability measure
$\mu_{\psi^{(R)}} = l r$ (i.e 
$\forall w \in \Omega^R$, $\mu_{\psi^{(R)}}(w)=l_w r_w$).
From this, one can compute the probability of spike blocks of arbitrary length
by the Chapman-Kolmogorov formula:
$$\mu_{\psi^{(R)}}(\cyl{\omega}{s}{t+R})=\mu_{\psi^{(R)}}(w(s))\prod_{n=s}^{t-1} \cL^{(R)}_{w(n)w(n+1)},$$
with $w(n) \sim \seq{\omega}{n}{n+R}$.

Then, one can check that $\mu_{\psi^{(R)}}$ is a Gibbs distribution \cite{keller:98} 
(see def. \ref{dGibbs}). Moreover, it is also an equilibrium state in the sense
of (\ref{VarPrinc}). As a result $0=h(\mu_{\psi^{(R)}})+\mu_{\psi^{(R)}}(\psi^{(R)})$,
therefore
$$h(\mu_{\psi^{(R)}})=-\mu_{\psi^{(R)}}(\psi^{(R)}).$$

\ssu{Convergence of the approximation.}\label{SRange}

Let us now discuss how well the range-$R$ potential (\ref{psi_R}) approximates
the infinite range potential (\ref{psi}). By definition of a range-$R+1$ potential
$\psi^{(R)}(\omega)=\psi^{(R)}(\seq{\omega}{-R}{0})$ so that we can compare
$\psi$ and $\psi^{(R)}$. 
One has:
$$\|\psi - \psi^{(R)} \|_\infty 
\leq
\sup\left\{ \, | \,\psi(\omega)-\psi(\omega')\, | \,: \omega, \omega' \in X, \omega(t)=\omega'(t), \, \forall t \in \left\{-R, \dots, 0 \right\}\right\} \deq var_R(\psi),$$
so that, from Theorem (\ref{ThUnicity}),
\beq\label{BoundsPsiR}
\|\psi - \psi^{(R)} \|_\infty \leq K' \gamma^R,
\eeq
where $K'$ is given by (\ref{Kp}). Therefore,
 $\psi^{(R)}$ approaches $\psi$ exponentially fast, as $R$ growths,
with a rate $\gamma$. This is in fact a classical result in ergodic theory:
regular potential are approximated by finite-range potential in the sup norm where
$\|\psi - \psi^{(R)} \|_\infty \leq C \Theta^R$, for some $0 < \Theta < 1$ (see
def. \ref{dX}). Here it is natural to take $\Theta=\gamma$.

The implications on statistics is related to the Kullback-Leibler divergence
(\ref{dKL}). Indeed, since $\mu_{\psi}$  and $\mu_{\psi^{(R)}}$ are
Gibbs distributions for the right shift $T$ we may use (\ref{dKLpres}), giving:
%
$$d\left(\mu_{\psi^{(R)}},\mu_\psi \right) = P(\psi) - \mu_{\psi^{(R)}}(\psi) - h(\mu_{\psi^{(R)}})
=\mu_{\psi^{(R)}}(\psi^{(R)}-\psi),
$$
%
where we used $P(\psi)=0$ (normalization of $\psi$) and $h(\mu_{\psi^{(R)}})=- \mu_{\psi^{(R)}}(\psi^{(R)})$
(see e.g. \cite{chazottes:99} for a more general proof).
Therefore,
\beq\label{ConvdKL}
d\left(\mu_{\psi^{(R)}},\mu_\psi \right) < K' \gamma^R.
\eeq

Therefore, the Kullback-Leibler divergence between the two measures
$\mu_\psi,\mu_{\psi^{(R)}}$, decays exponentially fast with a decay rate $\gamma$.

A practical consequence of this result is that it might be sufficient, for practical purposes, to approximate
$\psi$ with a potential of range:
\beq\label{RApprox}
R \sim -\frac{\log K'}{\log \gamma}.
\eeq
Note however that the constant $K'$ depend on several parameters. Especially,
it diverges when $\sigma_B \to 0$ (or $\gamma \to 1$). As a consequence,
depending on these parameters, the effective range can be quite large.

\su{Raster plots statistics}\label{SSpikeStat}

As discussed in the introduction, the neuroscience community is confronted to the delicate
problem of characterizing statistical properties of
raster plots from finite time spike trains and/or from finite number
of experiments. This requires an a priori guess
for the probability of raster plots, what we call 
a \textit{statistical model}. These models
can be extrapolated from heuristic arguments
or from principles such as the Jaynes
argument from statistical physics \cite{jaynes:57} (see
section \ref{SJaynes}).
In this section, we show that Markovian approximations
introduced in the previous section constitute such statistical models,
from which  classical
statistical indicators used by  the neuroscience community 
can be explicitly computed in the case of model (\ref{DNN}).

\ssu{Two representations of the potential $\psi^{(R)}$.}

\sssu{spike-block representation.}

The potential $\psi^{(R)}$ is a function of $W \sim \seq{\omega}{-R}{0}$.
Therefore, it 
takes only $L=2^{N(R+1)}$ values, \textit{explicitly given by (\ref{psi_R})}.
To each possible block $\seq{\omega}{-R}{0}$
we associate a word $W_n, n=1 \dots L$.
Call $\chi_n(W)$ the characteristic function,
equal to $1$, if $W=W_n$, and $0$ otherwise.
Then:
\beq\label{psi_alpha}
\psi^{(R)}(W)=\sum_{n=1}^{L=2^{N(R+1)}}\alpha_n\chi_n(W) \equiv \psi^{(R)}_\balpha(W), 
\eeq 
where $\alpha_n=\psi^{(R)}(W_n)$.
This decomposition of the potential is called the
\textit{spike-block representation} of $\psi^{(R)}$ and the index
$\balpha$ in (\ref{psi_alpha}) (which is the vector $\left(\alpha_n \right)_{n=1}^L$)
  makes this  representation explicit. 
Note that $\balpha$ depend (analytically) on the model-parameters 
$W_{ij}, i,j =1 \dots N, \, I_i, i=1 \dots N, \gamma, \theta$.

\sssu{Interpretation}

This representation is quite natural since $e^{\alpha_n}$ is nothing but
the probability $P(\omega(0) \, | \, \seq{\omega}{-R}{-1})$ 
with $W_n \sim \seq{\omega}{-R}{0}$, namely a matrix element of the
transition matrix (\ref{cL}). This function corresponds to the so-called
``conditional intensity'' introduced in neuroscience data analysis by researchers
like C. Pouzat and A. Chaffiol \cite{pouzat-chaffiol:09} and the exponential distribution
introduced by these authors is actually our Gibbs distribution. Here, we are able to compute
explicitly this distribution, because we are dealing with a model, while Pouzat
and Chaffiol are coping with real world data.
 
Fixing an history $\seq{\omega}{-R}{-1}$ the sum of  $e^{\alpha_n}$'s,
over all blocks $W_n$ having an history $\seq{\omega}{-R}{-1}$ and such that
$\omega_i(0)=1$, is the probability   
that neuron $i$ fires given the history $\seq{\omega}{-R}{-1}$.
More generally, the product of the transition matrix $\cL$ elements (\ref{cL}) provides the probability
of a certain sequence of spikes (``response'') given a certain history.
If one focuses on the response $\cR$ of a subset  of neurons in the network,
to spikes emitted by an another subset of neurons  in the network-corresponding
to a given history and considered as a stimulus $\cS$- the matrix $\cL$ in the
$\alpha$-representation allows the computation of  the probability $P(\cR \, | \, \cS)$. Then, by Bayesian
inference, and \textit{since the probability $P(\cS)$ of the stimulus is known}
(it is given by the invariant measure of the Markov chain), one infers $P(\cS \, | \, \cR)$.
This provides one way of characterizing the ``neural code'', in the sense
of \cite{rieke-etal:96}, at the level of networks of neurons, where stimuli
are spike trains.

\ssu{The spikes-uplets representation.}\label{lambda_rep} 

Though natural the $\alpha$-representation is not the most commonly used.
Let us introduce another representation.
 
\sssu{Monomials.}

An \textit{order-$n$ monomial} is
a product $\omega_{i_1}(t_1) \dots \omega_{i_n}(t_n)$,
where $1 \leq i_1 \leq i_2 \leq \dots \leq i_n \leq N$ and
$-\infty < t_1 \leq t_2 \leq \dots \leq t_n < +\infty$, and where there is
no repeated pair of indexes $(i,t)$.
Since  $\omega_i(t)^k=\omega_i(t), \forall i=1 \dots N, t \in \setZ, k>1$
the last requirement avoids redundancies.
A polynomial is a linear combination of monomials.

The monomial $\omega_{i_1}(t_1) \dots \omega_{i_n}(t_n)$
takes values in $\lbrace 0,1 \rbrace$ and is $1$ if and  only
if each neuron $i_l$ fires at times $t_l$, $l=1 \dots n$.
On phenomenological grounds this corresponds
to a spike $n$-uplet $(i_1,t_1), \dots, (i_n,t_n)$ (neuron $i_1$
fires at time $t_1$, and neuron $i_2$ fires at time $t_2$,
$\dots$).

\sssu{Spikes-uplets expansion of $\psi$.} 

Returning to the spike-block representation (\ref{psi_alpha}),
 the characteristic
function of the word $W_n$, $\chi_n(W)$, reads:
$$\chi_n(W)=Q_n(W)R_n(W),$$
with 
$$Q_n(W)=
\prod_{
(i,t),\omega_i(t)=1
}
\omega'_i(t),$$
where $W$ represents a spike block $\seq{\omega'}{-R}{0}$,
while the product holds over all pairs $(i,t)$, $1 \leq i \leq N$, $-R \leq t \leq 0$
such that $\omega_i(t)=1$, in the word $W_n \sim \seq{\omega}{-R}{0}$. Likewise,
%
$$R_n(W)=\prod_{(j,s), \omega_j(s)=0}(1-\omega'_j(s))
=\sum_{m=1}^{k_n}(-1)^m R_{n,m}(W),$$
%
where $R_{n,m}$ are monomials of order $\leq R$.
Since $\omega_i(t)^k=\omega_i(t), \, k \geq 1$,
 $Q_n$, $R_n$ are monomials of order $\leq R$
and all $\chi_n(W)$'s are polynomials of order $\leq R$.

For $N,R \in \setZ$, we note $\cP(N,R)$
the set of non repeated pairs of integers $(i,n)$ with
$i \in \left\{1,\dots,N \right\}$ and 
$t  \in \left\{-R,\dots,0\right\}$. We have just proved that
$\psi$ is approximated by a range-$R$  polynomial expansion of the form:
\beq\label{lambda_expansion}
\psi_\blambda^{(R)}(W)=
\sum_{l=0}^R
\sum_{
\tiny{
(i_1,t_1), \dots,(i_l,t_l) \in \cP(N,R),\\
}
}
\lambda^{(l)}_{i_1,t_1,\dots,i_l,t_l} 
\omega_{i_1}(t_1) \dots\omega_{i_l}(t_l),
\eeq
where $W \sim \seq{\omega}{0}{R}$. This is
called the ``spike-uplets representation''. It is obviously equivalent
to the spike-block representation $\psi^{(R)}_\balpha$
(the $\lambda_i$'s are linear combinations of the
$\alpha_l$'s) but this expansion is more convenient to discuss the link
between our results and the standard approaches used in the neuroscience
community. Note that a spike-block contains $0$'s and $1$'s (it tells us
which neurons are firing and which neurons are not firing) while
a spike-uplet takes only into account firing neurons. As a
consequence there is some redundancy in the spike-block representation
that can be removed in the spike-uplets representation (see details below). 

\sssu{Interpretation.}
Since the analytic function $\log(\pi(x))$ has a series
expansion for $x \in \setR$, setting 
$x=\frac{\theta-C_i(\seq{\omega}{-R}{-1})}{\sigma_i(\seq{\omega}{-R}{-1})}
$ in (\ref{psi_R}) $x^n$ is a sum of terms
$\omega_{i_1}(t_1) \dots \omega_{i_n}(t_n)$ and using the series 
expansion\footnote{Since $\omega_i(t)^k=\omega_i(t), \forall i=1 \dots N, t \in \setZ, k>1$,
 the terms of the series can be grouped
together giving rise to a finite sum
of monomials. 
}
one can compute explicitly the coefficients of the spike-uplets expansion.
Due to stationarity (see details below) one can only consider spike-uplets of the form 
$\omega_i(0)\omega_{j_1}(t_1) \dots \omega_{j_l}(t_l)$, with $t_1, \dots, t_l <0$.
They are 
combinations of terms proportional to $W_{i j_1}W_{i j_2} \dots W_{i j_l}\gamma^{-(t_1+t_2+ \dots t_l)}$
which have a nice interpretation. 
The sum of these terms corresponds  to 
the cumulative effect of 
spikes emitted by neurons $j_1, \dots, j_l$ at times $t_1, \dots, t_l$ in the past, on neuron $i$ at time $0$. 
Actually, these terms are related to a linear response theory as 
 developed in a different context in \cite{cessac-sepulchre:04,cessac-sepulchre:06}.

\ssu{Statistical models.}

\sssu{The topological pressure as a cumulant generator.}

Let us return to the  transition matrix $\cL^{(R)}$ and related
 topological pressure $P(\psi^{(R)})$ introduced in sections \ref{STrans},\ref{SPR}. 
These quantities depend on $\balpha$ or $\blambda$ according to the representation
(which is nothing but a change of variables).

The pressure is differentiable with respect to $\blambda$ and one has
\footnote{Here the pressure $P(\psi^{(R)})$ is 
 considered as a function of the $\blambda$'s, where these parameters are
arbitrary. As a consequence the corresponding potential $\psi_{\blambda}^{(R)}$
is no longer normalized. In this case the topological pressure is  defined as the logarithm of the maximal eigenvalue of the matrix ${\mathcal L}^{(R)}$. The fact that a potential of the form (\ref{lambda_expansion})
are in general \textit{not normalized} has deep practical consequences widely discussed in the paper \cite{vasquez-cessac-etal:10}.}:
%
$$\frac{\partial P(\psi^{(R)})}
{\partial  \lambda^{(l)}_{(i_1,t_1),\dots,(i_l,t_l)}} 
=\mu_{\psi^{(R)}} \left[\omega_{i_1}(t_1) \dots \omega_{i_l}(t_l) \right],
$$
%
Therefore, the derivation of the pressure with respect to the quantity
$\lambda^{(l)}_{(i_1,t_1),\dots,(i_l,t_l)}$ provides the $\mu_{\psi^{(R)}}$-probability
of the spike $n$-uplet $\omega_{i_1}(t_1) \dots \omega_{i_l}(t_l)$.  

In particular:
\beq\label{dPlambda}
\frac{\partial P(\psi^{(R)})}
{\partial \lambda^{(1)}_{i_1,t_1}} =
\mu_{\psi^{(R)}}\left[\omega_{i_1}(t_1)\right],
\eeq
the firing rate of neuron $i_1$ at time $t_1$.
Since dynamics is stationary this quantity does not depend
on $t$ (see eq. (\ref{defrate})). As a consequence all
terms $\lambda_{i,t}$, $i$ fixed and $t \in \left\{-R,\dots,0\right\}$
play the same role and we can simplify the potential 
(\ref{lambda_expansion}) in keeping, as first order
terms, the monomials of form $\lambda^{(1)}_i \omega_i(0), \, i=1 \dots N$.

In the same way:
%
$$\frac{\partial P(\psi^{(R)})}
{\partial  \lambda^{(2)}_{(i_1,t_1),(i_2,t_2)}} 
=\mu_{\psi^{(R)}} \left[\omega_{i_1}(t_1)\omega_{i_2}(t_2) \right],
$$
%
From stationarity
it follows that this quantity depends only on $t_2 -t_1$. So,
there are redundant terms in the expansion (\ref{lambda_expansion})
and we may write the part of the expansion corresponding to pairs of spikes
as $\sum_{\tau=-R/2}^{R/2} \sum_{i_1,i_2=1}^N \lambda_{i_1,i_2,\tau}^{(2)}\omega_{i_1}(0)\omega_{i_2}(\tau).$

Higher order
redundant terms can be removed as well, taking into 
account the stationarity of the process.
As a consequence we may write the  spike-uplets expansion of $\psi$ in
the form:
\beq\label{gener_psi_R}
\psi_\blambda^{(R)}(W)=\sum_{l=0}^K \lambda_l \phi_l(W),
\eeq
where 
$l$ enumerates all non redundant monomials $\phi_l$ of order
$\leq R$,
including the constant monomial $\phi_0(W)=1$.

\sssu{Statistics of raster plots from Jaynes formalism}\label{SJaynes}

We would like now to relate the present analysis to a standard problem in spike train analysis. Assume that we have generated a (finite)
raster plot $\omega_{exp}$
from the dynamical system (\ref{DNN}) and that we want to 
recover the probability distribution $\mu_\psi$ from this raster plot,
without any other information.
A usual approach consist of computing the average value of
some prescribed
spike-uplets, and infer the corresponding probability distribution
from a variational principle introduced by Jaynes \cite{jaynes:57}.
The Jaynes approach has been used by several authors in the field of
experimental spike train analysis \cite{schneidman-etal:06,tkacik-etal:06,marre-boustani-etal:09}.

 So, let us assume that, from the raster plot $\omega_{exp}$,
we have computed, by time average\footnote{This argument extends
to the case where several empirical raster plots have been generated.
Then, average values are obtained by combinations of time average
and sample average.}, 
 the average value $C_l$ of an a priori fixed set of 
monomials $\phi_l$, $l=1 \dots M$.
To find a probability distribution $\mu_{\psi_{test}}$ which matches these
average values, without making additional assumptions,
one maximizes the statistical entropy under the constraints $\mu_{\psi_{test}}(\phi_l)=C_l$, $l=1 \dots M$.
In the context of  thermodynamic formalism  \textit{this amounts
to finding a set of parameters $\lambda_l$ satisfying the variational equation (\ref{VarPrinc})
for a finite range potential} $\psi_{test}=\sum_{l=1}^M \lambda_l \phi_l$.
The $\lambda_l$'s 
are adjustable Lagrange multipliers,
which have to be tuned using (\ref{dPlambda}), so that the 
 average of $\phi_l$ with respect to $\mu_{\psi_{test}}$  is equal to $C_l$

Let us give two classical examples.

\sssu{Homogeneous Bernoulli statistics.}  
The simplest example consists of only measuring, thus, constraining the value of firing 
rates.
This amount to considering a range-$1$ potential
\cite{rieke-etal:96,dayan-abbott:01,gerstner-kistler:02b}.
:
$$\psi_{test}(W)=  
\sum_{i=1}^N 
\lambda^{(1)}_{i} 
\omega_{i}(0),$$
corresponding to a  probability:
$$
\mu_{\psi_{test}}(\omega(0))=
\prod_{i=1}^N \frac{e^{\lambda^{(1)}_{i}\omega_{i}(0)}}{1+e^{\lambda^{(1)}_{i}}}.$$
Therefore, taking only the first order monomials
allows one to select a probability distribution under which
neurons fire independently with
a time-independent rate: 
$$r_i=
\frac{e^{\lambda^{(1)}_{i}}}{1+e^{\lambda^{(1)}_{i}}}.$$
\sssu{Pairwise interactions.}
  
This statistical model  has been introduced, in the
context of spike trains statistics, by Schneidman et al \cite{schneidman-etal:06}.
Here,  still $R=1$ and $\psi_{test}$ does not depend on $t$ but pairs
of spikes occurring at the same time are considered. 
Then, 

$$ 
\psi_{test}(W) =
\sum_{i=1}^N 
\lambda^{(1)}_{i} 
\omega_{i}(0)
+
\sum_{1 \leq i_1<i_2 \leq N} 
\lambda^{(2)}_{i_1,i_2} 
\omega_{i_1}(0)\omega_{i_2}(0).$$

Here, the related probability measure does not factorize any more
but all information about spike train statistics is contained
in the first and second order spike-uplets.

\sssu{Choosing an a priori set of monomials.}

More general potentials can be considered as well \cite{marre-boustani-etal:09}. In view of the present analysis, fixing an a priori
set of observables, often fixed from  a priori hypotheses
on the relative role of spike-uplets (e.g. rate versus synchronization), amounts to fixing a test potential $\psi_{test}$. 
Therefore, there are as many models as possible choices of observables.
How to discriminate them ?

The probability distribution $\mu_{\psi_{test}}$ is the Gibbs distribution
for the potential $\psi_{test}$. It provides an approximation of
the invariant measure $\mu_{\psi}$ in two senses. First,
$\psi_{test}$ contains only some terms in the polynomial expansion
of a finite range potential $\psi^{(R)}$, which are fixed from the a priori choice of observables.
Second,  $\psi^{(R)}$ is  a finite-range approximation of
the exact, infinite-range, potential $\psi$. 
The ``error'' is measured by the Kullback-Leibler divergence (\ref{dKL}):
 $$d\left(\mu_{\psi_{test}},\mu_\psi \right) =
\mu_{\psi_{test}}(\psi_{test}-\psi).
$$
It is upper-bounded by $C\gamma^{NR}$. This fixes, for model (\ref{DNN}), an estimate for
the value of $R$ given by (\ref{RApprox}) .

Nevertheless, the number of terms increases exponentially
with $R$ and $N$ and therefore, especially if $\gamma$ is close
to $1$ there is an overwhelming number of monomials. 
 Now, it might be,  that some terms $\phi_l$ are less
important than others: the corresponding coefficient $\lambda_l$
vanishes or is small compared to others terms.
As discussed in \cite{cessac-rostro-etal:09} and shortly commented
in section \ref{SPlast}, neural mechanisms such as plasticity
certainly reinforce some of these terms (especially rates and spike pairs). On a more abstract setting,
our analysis shows
that the Kullback-Leibler divergence (\ref{dKL}) gives
an indication of the distance between the probability 
reconstructed from Jaynes principle, with a ``guess'' potential, and the true probability $\mu_\psi$. This opens up a way to compare and select statistical models by minimizing the Kullback-Leibler divergence,
using eq. (\ref{dKL}). This aspect and its numerical implementation
are discussed in \cite{vasquez-cessac-etal:10}.

\su{Discussion and conclusion}

In this paper we have addressed the question of characterizing the spike train statistics of a network of LIF neurons with noise, in the stationary case,
 with two aims. Firstly, to obtain analytic and rigorous results allowing the characterization of the process of spike generations. For this, we have used the realm of ergodic theory and thermodynamic formalism, which looks well adapted for this purpose. We have obtained unexpected results, especially the fact that, even  so simple models of neural networks have, strictly speaking, an unbounded memory rendering spike train statistics non-Markovian. The common wisdom in the field of neural networks dynamics suggests, however, that there is a characteristic time scale after which the system essentially looses its memory. Here, this time scale is controlled by $\gamma$, the leak rate, closely related
to synaptic response time. 

The second goal was to make a connection from this mathematical analysis toward the empirical methods used in neuroscience community for the analysis of spike trains. Here, we have shown that the Jaynes method, based on an a priori choice of a ``guess'' potential, with finite range, amounts to approximate the exact probability distribution by the Gibbs distribution of a Markov chain \cite{csiszar:84}. The degree of approximation
can be controlled by the Kullback-Leibler divergence which can computed using a classical result in the thermodynamic formalism. This analysis opens up the possibility of developing efficients algorithms
to estimate at best the statistic of spike trains from experimental data, using several guess potential and selecting the one which minimizes the KL divergence \cite{vasquez-cessac-etal:10}.\\
 
Clearly, this work is just a beginning, since, especially, it deals with a rather simple model.
Let us now briefly comment several possible extensions.

\ssu{Conductance based Integrate-and-Fire neurons.}

A natural extension of the present works concerns the so-called Generalized Integrate-and-Fire
models \cite{rudolph-destexhe:06} , which are closer to
 biology \cite{jolivet-lewis-etal:04,jolivet-rauch-etal:06}.
The occurrence of a post-synaptic potential on synapse $j$, at time $\tjn$,
results in a change of membrane potential. In conductance based models 
this change is integrated in the adaptation of conductances. It has been
shown in  \cite{cessac-vieville:08} that, under natural assumptions on
spike-time precision that the continuous-time evolution of these equations
reduces to the discrete time dynamics:

$$V_i(t+1)= \gamma_i(t,\tot)\left[1-Z(V_i(t))\right]V_i(t)+\Jito, \quad i=1 \dots N,$$

\nid where:

$$\gamma_i(t,\tot) \deq e^{-\int_{t}^{t+1} \, g_i(s,\tot) \, ds} < 1,$$

\nid is the integral of the  conductance $g_i(s,\tot)$ over the time interval $[t,t+1[$.
Conductances depend on the past spikes via the relation:

$$g_{ij}(t,\tot)=  G_{ij} \sum_{n=1}^{M_j(\seq{\omega}{-\infty}{t})} \alpha_j(t-\tjn).$$ 

\nid In this equation, $M_j(\seq{\omega}{-\infty}{t})$ is the number
of times neuron $j$ has fired at time $t$ (it can be infinite).
$\alpha$ is the synaptic profile (it decays exponentially fast)
and $\tjn$ is the time of occurrence of the $n$-th spike in the raster
$\omega$.
  $G_{ij}$ is a positive constant proportional to
 the synaptic efficacy
$$\left\{
\baR{ccc}
 W_{ij}=E^+G_{ij} \quad &\mbox{if}&  \quad j \in \cE, \\
 W_{ij}=E^-G_{ij} \quad &\mbox{if}&  \quad j \in \cI, 
\eaR
\right. 
$$
where $E_L,E^{+},E^{-}$ are respectively the Nernst potentials for the leak,
the excitatory (set $\cE$) and the inhibitory synapses (set $\cI$).

The term,
$$\Jito=\int_{t}^{t+1} i_i(s,\tot) \, \nu_i(s,t+1,\tot) \, ds,$$
\nid is the corresponding integrated synaptic current with:
$$i_i(t,\tot)=\frac{E_L}{\tau_L} 
+  E^{+} \, \sum_{j \in \cE} g_{ij}(t,\tot)
+  E^{-} \, \sum_{j \in \cI} g_{ij}(t,\tot)
+i^{(ext)}_i(t),
$$
$$\nu_i(s,t+1,\tot)=e^{-\int_{s}^{t+1} \, g_i(s',\tot) \, ds'}.$$.

The difficulty here is that the coefficient $\gamma_i(t,\tot)$, which is the analog
of $\gamma$ in eq. (\ref{DNN}) depends on the whole past. This introduces
another non-Markovian effect in the dynamics. In this case the computation
of the potential corresponding to (\ref{psi}) is clearly more complex. 
This case is under current investigations.

\ssu{Non stationarity.}

One weakness of the present work is that it only considers
stationary dynamics, where e.g. the external current $I_i$
is independent of time. Besides, we have taken the limit
$s \to -\infty$ in section \ref{Sdef} to remove the dependence in the
initial condition $V(s)$. However, real neural systems are submitted
to non static stimuli, and transients play a crucial role.
To extend the present analysis to these case one needs the
proper mathematical framework. The non stationarity requires
to handle time dependent Gibbs measures. In the realm of
ergodic theory applied to non equilibrium statistical
physics, Ruelle has introduced the notion of time-dependent
SRB measure \cite{ruelle:99}. A similar approach could be used
here, at least formally.

Handling the transients is an even more tricky question. The main
difficulty is to propose a probability distribution for the initial
condition $V(s)$. From the dynamical systems point of view it is natural to take e.g. Lebesgue, and extensions toward
this case are under current investigations. But if one wants to make serious
extrapolations of mathematical results towards neuroscience one has
to ask why the ``initial state'' of a neural network, namely the
state in which the neural network is as the experiment starts, should be uniform in the phase space (or Gaussian or 
whatsoever), as soon as this initial state is the result
of a previous (phylogenetic and ontogenetic) evolution ?

\ssu{Synaptic plasticity.}\label{SPlast}

In neural networks, synaptic weights are not fixed, as in (\ref{DNN}), but they evolve with the activity of the pre- and post-synaptic
neuron (synaptic plasticity). This means that synaptic weights evolve
according to spike train statistics, while spike train statistics is
constrained by synaptic weights. This interwoven evolution has been considered in \cite{cessac-rostro-etal:09} under the assumption that spike-train statistics is characterized by a Gibbs distribution. Actually, the present work confirms this hypothesis in the case of LIF models. The main conclusion of \cite{cessac-rostro-etal:09} is that synaptic mechanism occurring on a time scale which is slow compared to neural dynamics are associated with a variational principle. There is a function, closely related to the topological pressure, which decreases when the synaptic adaptation process takes place. Moreover, the synaptic adaptation has the effect of reinforcing specific terms in the potential, directly related to the form of the synaptic plasticity mechanism.
The interest of this result is that it provides an a priori guess
of the relevant terms in the potential expansion. A contrario, it allows one to constrain the spike train statistics of a LIF model, using synaptic plasticity with an appropriate rule which can be determined
from the form of the expected potential.

Finally, an interesting issue fitting together with the discussion of non stationarity and synaptic plasticity,
is to analyse spike frequency adaptation in this context
\cite{crook-ermentrout-etal:98,ermentrout-pascal-etal:01,benda-hertz:03,benda-longtin-etal:05}

\begin{acknowledgements}
I am grateful to reviewers for a careful reading of the manuscript, constructive criticism
and helpful comments.
I would like to thank J.R. Chazottes, O. Faugeras, B. Fernandez, F. Grammont, J. Touboul, J.C. Vasquez, T. Vi\'eville for helpfull discussions. This work has been partially supported by the INRIA ARC grant MACACC.

\end{acknowledgements}

\bibliographystyle{abbrv}

\bibliography{biblio,odyssee}

\begin{thebibliography}{10}

\bibitem{adrian-zotterman:1926}
E.~Adrian and Y.~Zotterman.
\newblock The impulses produced by sensory nerve endings: Part ii: The response
  of a single end organ.
\newblock {\em J Physiol (Lond.)}, 61:151--71, 1926.

\bibitem{benda-hertz:03}
J.~Benda and A.~Herz.
\newblock {A universal model for spike-frequency adaptation}.
\newblock {\em {NEURAL COMPUTATION}}, {15}({11}):{2523--2564}, {NOV} {2003}.

\bibitem{benda-longtin-etal:05}
J.~Benda, A.~Longtin, and L.~Maler.
\newblock Spike-frequency adaptation separates transient communication signals
  from background oscillations.
\newblock {\em J. Neurosci.}, 25:2312 -- 2321, 2005.

\bibitem{bowen:98}
R.~Bowen.
\newblock {\em Equilibrium states and the ergodic theory of Anosov
  diffeomorphisms. Second revised version.}, volume 470 of {\em Lect. Notes.in
  Math.}
\newblock Springer-Verlag, 2008.

\bibitem{bressaud-fernandez-etal:99}
X.~Bressaud, R.~Fernandez, and A.~Galves.
\newblock Decay of correlations for non h\"olderian dynamics. a coupling
  approach.
\newblock {\em Electronic Journal of Probabilities}, 4(3):1--19, 1999.

\bibitem{brunel-hakim:99}
N.~Brunel and V.~Hakim.
\newblock Fast global oscillations in networks of integrate-and-fire neurons
  with low firing rates.
\newblock {\em Neural Computation}, 11:1621--1671, 1999.

\bibitem{cessac:08}
B.~Cessac.
\newblock A discrete time neural network model with spiking neurons. i.
  rigorous results on the spontaneous dynamics.
\newblock {\em J. Math. Biol.}, 56(3):311--345, 2008.

\bibitem{cessac-rostro-etal:09}
B.~Cessac, H.~Rostro-Gonzalez, J.~Vasquez, and T.~Vi\'eville.
\newblock How gibbs distribution may naturally arise from synaptic adaptation
  mechanisms: a model based argumentation.
\newblock {\em J. Stat. Phys}, 136(3):565--602, August 2009.

\bibitem{cessac-sepulchre:04}
B.~Cessac and J.~Sepulchre.
\newblock Stable resonances and signal propagation in a chaotic network of
  coupled units.
\newblock {\em Phys. Rev. E}, 70(056111), 2004.

\bibitem{cessac-sepulchre:06}
B.~Cessac and J.~Sepulchre.
\newblock Transmitting a signal by amplitude modulation in a chaotic network.
\newblock {\em Chaos}, 16(013104), 2006.

\bibitem{cessac-vieville:08}
B.~Cessac and T.~Vi{\'e}ville.
\newblock On dynamics of integrate-and-fire neural networks with adaptive
  conductances.
\newblock {\em Frontiers in neuroscience}, 2(2), jul 2008.

\bibitem{chazottes:99}
J.~Chazottes.
\newblock {\em Entropie Relative, Dynamique Symbolique et Turbulence}.
\newblock PhD thesis, Universit\'e de Provence - Aix Marseille I, 1999.

\bibitem{chazottes-keller:09}
J.~Chazottes and G.~Keller.
\newblock {\em Pressure and Equilibrium States in Ergodic Theory}, chapter
  Ergodic Theory.
\newblock Encyclopedia of Complexity and System Science, Springer, 2009.
\newblock to appear.

\bibitem{coelho-quas:98}
Z.~Coelho and A.~Quas.
\newblock criteria for $d$-continuity.
\newblock {\em Transactions of the American Mathematical Society},
  350(8):3257--3268, 1998.

\bibitem{crook-ermentrout-etal:98}
S.~Crook, G.~Ermentrout, and J.~Bower.
\newblock {Spike frequency adaptation affects the synchronization properties of
  networks of cortical oscillators}.
\newblock {\em {NEURAL COMPUTATION}}, {10}({4}):{837--854}, {MAY 15} {1998}.

\bibitem{csiszar:84}
I.~Csiszar.
\newblock Sanov property, generalized $i$-projection and a conditional limit
  theorem.
\newblock {\em Ann. Prob}, 12(3):768--793, 1984.

\bibitem{dayan-abbott:01}
P.~Dayan and L.~F. Abbott.
\newblock {\em Theoretical Neuroscience : Computational and Mathematical
  Modeling of Neural Systems}.
\newblock MIT Press, 2001.

\bibitem{ermentrout-pascal-etal:01}
B.~Ermentrout, M.~Pascal, and B.~Gutkin.
\newblock {The effects of spike frequency adaptation and negative feedback on
  the synchronization of neural oscillators}.
\newblock {\em {NEURAL COMPUTATION}}, {13}({6}):{1285--1310}, {JUN} {2001}.

\bibitem{gantmacher:66}
F.~R. Gantmacher.
\newblock {\em the theory of matrices}.
\newblock AMS Chelsea Publishing, Providence, RI, 1998.

\bibitem{gerstner-kistler:02b}
W.~Gerstner and W.~Kistler.
\newblock {\em Spiking Neuron Models}.
\newblock Cambridge University Press, 2002.

\bibitem{gerstner-kistler:02}
W.~Gerstner and W.~M. Kistler.
\newblock Mathematical formulations of hebbian learning.
\newblock {\em Biological Cybernetics}, 87:404--415, 2002.

\bibitem{grammont-riehle:99}
F.~Grammont and A.~Riehle.
\newblock Precise spike synchronization in monkey motor cortex involved in
  preparation for movement.
\newblock {\em Exp. Brain Res.}, 128:118--122, 1999.

\bibitem{grammont-riehle:03}
F.~Grammont and A.~Riehle.
\newblock Spike synchronization and firing rate in a population of motor
  cortical neurons in relation to movement direction and reaction time.
\newblock {\em Biol Cybern}, 88:360--373, 2003.

\bibitem{izhikevich:04}
E.~Izhikevich.
\newblock Which model to use for cortical spiking neurons?
\newblock {\em IEEE Trans Neural Netw}, 15(5):1063--1070, September 2004.

\bibitem{jaynes:57}
E.~Jaynes.
\newblock Information theory and statistical mechanics.
\newblock {\em Phys. Rev.}, 106(620), 1957.

\bibitem{johansson-oberg:2002}
A.~Johansson and A.~Oberg.
\newblock Square summability of variations of $g$-functions and uniqueness of
  $g$-measure.
\newblock {\em Math. res. Lett.}, 10:587--601, 2003.

\bibitem{jolivet-lewis-etal:04}
R.~Jolivet, T.~Lewis, and W.~Gerstner.
\newblock Generalized integrate-and-fire models of neuronal activity
  approximate spike trains of a detailed model to a high degree of accuracy.
\newblock {\em Journal of Neurophysiology}, 92:959--976, 2004.

\bibitem{jolivet-rauch-etal:06}
R.~Jolivet, A.~Rauch, H.-R. Lüscher, and W.~Gerstner.
\newblock Predicting spike timing of neocortical pyramidal neurons by simple
  threshold models.
\newblock {\em Journal of Computational Neuroscience}, 21:35--49, 2006.

\bibitem{keane:72}
M.~Keane.
\newblock Strongly mixing g-measures.
\newblock {\em Invent. Math.}, 16:309--324, 1972.

\bibitem{keller:98}
G.~Keller.
\newblock {\em Equilibrium States in Ergodic Theory}.
\newblock Cambridge University Press, 1998.

\bibitem{kirst-geisel-etal:09}
C.~Kirst, T.~Geisel, and M.~Timme.
\newblock Sequential desynchronization in networks of spiking neurons with
  partial reset.
\newblock {\em Phys. Rev. Lett.}, 102(068101), 2009.

\bibitem{kitchens:98}
B.~P. Kitchens.
\newblock {\em Symbolic Dynamics: One-sided, Two-sided and Countable State
  Markov Shifts}.
\newblock Springer-Verlag, 1998.

\bibitem{lapicque:07}
L.~Lapicque.
\newblock Recherches quantitatives sur l'excitation \'electrique des nerfs
  trait\'ee comme une polarisation.
\newblock {\em J. Physiol. Pathol. Gen.}, 9:620--635, 1907.

\bibitem{ledrappier:74}
F.~Ledrappier.
\newblock Principe variationnel et syst\`emes dynamiques symboliques.
\newblock {\em Z. Wahr. verw. Gebiete}, 30(185):185--202, 1974.

\bibitem{maillard:07}
G.~Maillard.
\newblock {\em Introduction to chains with complete connections}.
\newblock Ecole Federale Polytechnique de Lausanne, winter 2007.

\bibitem{marre-boustani-etal:09}
O.~Marre, S.~E. Boustani, Y.~Fr\'egnac, and A.~Destexhe.
\newblock Prediction of spatiotemporal patterns of neural activity from
  pairwise correlations.
\newblock {\em Phys. rev. Let.}, 102:138101, 2009.

\bibitem{nirenberg-latham:03}
S.~Nirenberg and P.~Latham.
\newblock Decoding neuronal spike trains: how important are correlations.
\newblock {\em Proceeding of the Natural Academy of Science},
  100(12):7348--7353, 2003.

\bibitem{onicescu-mihoc:35}
O.~Onicescu and G.~Mihoc.
\newblock Sur les cha\^ines statistiques.
\newblock {\em C. R. Acad. Sci. Paris}, 200, 1935.

\bibitem{perrinet:08}
L.~Perrinet.
\newblock Sparse spike coding : applications of neuroscience to the processing
  of natural images.
\newblock In S.~of~Photo-Optical Instrumentation~Engineers, editor, {\em
  Proceedings of SPIE, the International Society for Optical Engineering},
  number ISSN 0277-786X CODEN PSISDG, Bellingham, WA, ETATS-UNIS, 2008.

\bibitem{pouzat-chaffiol:09}
C.~Pouzat and A.~Chaffiol.
\newblock On goodness of fit tests for models of neuronal spike trains
  considered as counting processes.
\newblock {\em http://arxiv.org/abs/0909.2785v1}, 2009.

\bibitem{riehle-etal:00}
A.~Riehle, F.~Grammont, M.~Diesmann, and S.~Grün.
\newblock Dynamical changes and temporal precision of synchronized spiking
  activity in monkey motor cortex during movement preparation.
\newblock {\em J. Physiol (Paris)}, 94:569--582, 2000.

\bibitem{rieke-etal:96}
F.~Rieke, D.~Warland, R.~de~Ruyter~van Steveninck, and W.~Bialek.
\newblock {\em Spikes, Exploring the Neural Code}.
\newblock The M.I.T. Press, 1996.

\bibitem{roudy-nirenberg-etal:09}
Y.~Roudy, S.~Nirenberg, and P.~Latham.
\newblock Pairwise maximum entropy models for studying large biological
  systems: when they can work and when they can't.
\newblock {\em PLOS Computational Biology}, 5(5), 2009.

\bibitem{rudolph-destexhe:06}
M.~Rudolph and A.~Destexhe.
\newblock Analytical integrate and fire neuron models with conductance-based
  dynamics for event driven simulation strategies.
\newblock {\em Neural Computation}, 18:2146--2210, 2006.

\bibitem{ruelle:69}
D.~Ruelle.
\newblock {\em Statistical Mechanics: Rigorous results}.
\newblock Benjamin, New York, 1969.

\bibitem{ruelle:78}
D.~Ruelle.
\newblock {\em Thermodynamic formalism}.
\newblock Addison-Wesley,Reading, Massachusetts, 1978.

\bibitem{ruelle:99}
D.~Ruelle.
\newblock Smooth dynamics and new theoretical ideas in nonequilibrium
  statistical mechanics.
\newblock {\em J. Statist. Phys.}, 95:393--468, 1999.

\bibitem{rullen-thorpe:01}
R.~V. Rullen and S.~Thorpe.
\newblock Rate coding versus temporal order coding: What the retina ganglion
  cells tell the visual cortex.
\newblock {\em Neural Computing}, 13(6):1255--1283, 2001.

\bibitem{schneidman-etal:06}
E.~Schneidman, M.~Berry, R.~Segev, and W.~Bialek.
\newblock Weak pairwise correlations imply string correlated network states in
  a neural population.
\newblock {\em Nature}, 440:1007-- 1012, 2006.

\bibitem{seneta:06}
E.~Seneta.
\newblock {\em Non-negative Matrices and Markov Chains}.
\newblock Springer, 2006.

\bibitem{sinai:72}
Y.~Sinai.
\newblock Gibbs measures in ergodic theory.
\newblock {\em Russ. Math. Surveys}, 27(4):21--69, 1972.

\bibitem{soula-beslon-etal:06}
H.~Soula, G.~Beslon, and O.~Mazet.
\newblock Spontaneous dynamics of asymmetric random recurrent spiking neural
  networks.
\newblock {\em Neural Computation}, 18(1), 2006.

\bibitem{thorpe:90}
S.~Thorpe.
\newblock Spike arrival times: A highly efficient coding scheme for neural
  networks.
\newblock {\em Parallel processing in neural systems and computers}, pages
  91--94, 1990.

\bibitem{thorpe-fize-etal:96}
S.~Thorpe, D.~Fize, and C.~Marlot.
\newblock Speed of processing in the human visual system.
\newblock {\em Nature}, 381:520--522, 1996.

\bibitem{tkacik-etal:06}
G.~Tkacik, E.~Schneidman, M.~Berry, and W.~Bialek.
\newblock Ising models for networks of real neurons.
\newblock {\em arXiv}, q-bio/0611072, 2006.

\bibitem{touboul-faugeras:07}
J.~Touboul and O.~Faugeras.
\newblock The spikes trains probability distributions: a stochastic calculus
  approach.
\newblock {\em Journal of Physiology, Paris}, 101/1-3:78--98, dec 2007.

\bibitem{touboul-faugeras:09}
J.~Touboul and O.~Faugeras.
\newblock A markovian event-based framework for stochastic spiking neural
  networks.
\newblock Technical report, arXiv, 2009.
\newblock Submitted to Neural Computation.

\bibitem{toyoizumi-etal:07}
T.~Toyoizumi, J.-P. Pfister, K.~Aihara, and W.~Gerstner.
\newblock Optimality model of unsupervised spike-timing dependent plasticity:
  Synaptic memory and weight distribution.
\newblock {\em Neural Computation}, 19:639--671, 2007.

\bibitem{vasquez-cessac-etal:10}
J.~Vasquez, B.~Cessac, and T.~Vieville.
\newblock Entropy-based parametric estimation of spike train statistics.
\newblock {\em Journal of Computational Neuroscience}, 2010.
\newblock submitted.

\end{thebibliography}

\end{document}